\documentclass[a4wide,reqno,11pt]{amsart}
\pdfoutput=1

\usepackage[utf8]{inputenc}

\usepackage{amsmath,amsthm,amssymb,amsfonts,amscd,amsxtra}
\usepackage{a4wide} 
\usepackage[]{hyperref}
\usepackage{graphicx}
\usepackage[all]{xy}  
\usepackage{latexsym} 
\usepackage{caption}
\usepackage{soul}

\newtheorem{thm}{Theorem}[section]
\newtheorem{prop}[thm]{Proposition}
\newtheorem{lem}[thm]{Lemma}

\newtheorem*{claim*}{Claim}

\theoremstyle{definition}
\newtheorem{defn}{Definition}[section]
\newtheorem{example}{Example}

\newtheorem*{examples*}{Examples}
\newtheorem*{example*}{Example}

\theoremstyle{remark}
\newtheorem{rmk}[thm]{Remark}

\renewcommand{\natural}{\sharp}

\newcommand{\x}{\mathbf x}
\newcommand{\y}{\mathbf y}

\title{Lagrangian antisurgery}
\author{Luis Haug}
\address{Centre de Recherches Mathématiques, Université de Montréal,
    Pavillon André-Aisenstadt, 2920 Chemin de la tour, Montréal (QC),
    H3T 1J4, Canada}
\email{lhaug@posteo.net}
\date{\today}

\usepackage{color}

\begin{document}
\begin{abstract}
    We describe an operation which modifies a Lagrangian submanifold
    $L$ in a symplectic manifold $(M, \omega)$ such as to produce a
    new immersed Lagrangian submanifold $L'$, which as a smooth
    manifold is obtained by surgery along a framed sphere in
    $L$. Intuitively, this can be described as collapsing an isotropic
    disc with boundary on $L$ to a point. The inverse operation
    generalizes classical Lagrangian surgery. We also describe
    corresponding immersed Lagrangian cobordisms between $L$ and
    $L'$. After removal of their singular locus, we obtain examples of
    embedded Lagrangian cobordisms with precisely two ends. As an
    application, we use this construction to produce interesting
    examples of Lagrangian cobordisms between Clifford and Chekanov
    tori.
\end{abstract}
\maketitle

\section{Introduction}
\label{sec:introduction}
A fundamental question in symplectic geometry is what manifolds arise
as the Lagrangian submanifolds of a given symplectic manifold
$(M^{2n},\omega)$. This question has different flavours and levels of
difficulty, depending on whether one asks for embedded or immersed
Lagrangian submanifolds, and on whether one incorporates constraints
such as exactness or monotonicity.

A natural attempt to construct new Lagrangian submanifolds is to
modify given ones by some sort of surgery operation. There is one well
known construction which resolves the transverse double points of a
Lagrangian immersion $\iota: L \to M$ by replacing neighbourhoods of
them by copies of $D^1 \times S^{n-1}$. For example, if $L$ is
connected, oriented and immersed with a unique double point, then the
resulting Lagrangian $L'$ is embedded and diffeomorphic to the
connected sum $L \# (S^1 \times S^{n-1})$, provided that the surgery
can be performed compatibly with the orientation. This operation,
which we will refer to as \emph{Lagrangian 0-surgery}, is due to
Lalonde--Sikorav \cite{Lalonde-Sikorav91} for $n = 2$ and to
Polterovich \cite{Polterovich--Surgery-of-Lagrangians} for general
$n$.

\subsection*{Terminology and notation} In all of the following,
``Lagrangian submanifolds'' will generally be allowed to be immersed
with transverse double points. We will usually not make a notational
distinction between abstract smooth manifolds $L$ and their immersed
images in $M$; that is, whenever we have a Lagrangian immersion
$\iota: L \to M$, we will slightly abuse notation and denote its image
$\iota(L) \subset M$ also by $L$.

\subsection{Surgery of smooth manifolds}
\label{sec:surgery}
On the level of abstract smooth manifolds (i.e., not taking into
account the Lagrangian embedding), the passage from $L$ to $L'$ by
Lagrangian 0-surgery replaces an embedded copy of $S^0 \times D^n$ by
a copy of $D^1 \times S^{n-1}$. This is a special case of the
following more general operation originally due to Milnor
\cite{Milnor--Killing-htpy-groups61}: Whenever a smooth
$n$-dimensional manifold $L$ contains an embedding
$\varphi: S^k \times D^{n-k} \to L$, one can cut out
$\varphi(S^k \times D^{n-k})$ and replace it by a copy of
$D^{k+1} \times S^{n-k-1}$, such as to obtain a new manifold
\begin{equation*}
    L' = (L \smallsetminus \varphi(S^k \times D^{n-k}))
    \cup_{\varphi(S^k \times S^{n-k-1})} (D^{k+1} \times S^{n-k-1}).
\end{equation*}
This works because $\partial (S^k \times D^{n-k}) = S^k \times
S^{n-k-1} = \partial (D^{k+1} \times S^{n-k-1})$.  We say that the
manifold $L'$, which inherits a smooth structure from $L$ in a
canonical way, is obtained from $L$ by \emph{$k$-surgery} (a.k.a.\
\emph{surgery of index $k+1$}).

Surgery theory is closely connected to cobordism theory. The manifold
$L'$ resulting from $k$-surgery on a manifold $L$ is cobordant to $L$
via a cobordism
\begin{equation*}
    V = ([0,1] \times L) \cup_{\{1\} \times \varphi(S^k \times
        D^{n-k})} D^{k+1} \times D^{n-k},
\end{equation*}
i.e., a cobordism that arises from the cylinder $[0,1] \times L$ by
attaching a $(k+1)$-handle $D^{k+1} \times D^{n-k}$ along $\{1\}
\times \varphi(S^{k} \times D^{n-k})$. This cobordism is called the
\emph{trace} of the corresponding surgery.

\subsection{Lagrangian antisurgery}
\label{sec:lagr-surg-high}
Let now $L \subset M$ be a Lagrangian submanifold containing an
embedded copy of $S^k \times D^{n-k}$. It is natural to ask if the
manifold $L'$ obtained by $k$-surgery on $L$ can again be embedded or
immersed into $M$ as a Lagrangian submanifold. The answer to a
sufficiently strong version of this question is certainly negative:
For example, a closed orientable manifold $L$ that can be embedded in
$\mathbb C^n$ must have Euler characteristic $\chi(L) = 0$. However,
$k$-surgery changes the Euler characteristic according to
$\chi(L') = \chi(L) + (-1)^{k+1} + (-1)^{n-k-1}$, and hence does not
preserve its vanishing if $n$ is even. So in this case no result of a
single $k$-surgery on $L$ admits a Lagrangian embedding into
$\mathbb C^n$.

In this paper we will describe a construction which implements
$k$-surgery for Lagrangian submanifolds under certain conditions. Let
$L \subset M$ be a Lagrangian submanifold containing an embedding
$\varphi: S^k \times D^{n-k} \to L$ together with an \emph{isotropic
    surgery disc} $D$, that is, an embedded isotropic $(k+1)$-disc
$D \subset M$ intersecting $L$ cleanly along
$S = \varphi(S^k \times \{0\})$ and otherwise disjoint from $L$ (this
terminology is borrowed from \cite{DimRizell--LegAmbSurg}).

\begin{thm}
    \label{main-thm:antisurgery}
    The manifold $L'$ obtained by $k$-surgery on $L$ with respect to
    the embedding $\varphi: S^k \times D^{n-k} \to L$ admits a
    Lagrangian immersion $L' \to M$ whose image agrees with $L$
    outside of an arbitrarily small neighbourhood of $D$, and such
    that in this neighbourhood it has exactly one transverse double
    point. Moreover, there exists an immersed Lagrangian cobordism
    $ V: L' \leadsto L$ given by a Lagrangian immersion of the trace
    of the $k$-surgery into $T^*\mathbb R \times M$, whose singular
    locus is a 1-dimensional family of double points along the end
    corresponding to $L'$.
\end{thm}

The construction of $L'$ and $V: L' \leadsto L$, and hence the proof
of Theorem \ref{main-thm:antisurgery}, is the content of Section
\ref{sec:constr-cobord}. We refer to the operation that passes from
$L$ to $L'$ as \emph{Lagrangian k-antisurgery}. The idea behind the
terminology is that the operation \emph{creates} a double point, in
contrast to Lagrangian 0-surgery, which \emph{resolves} a double
point. To give a quick and intuitive description, one could say that
Lagrangian $k$-antisurgery modifies a Lagrangian $L$ by collapsing an
isotropic $(k+1)$-disc with boundary on $L$ to a point.

The local model for the immersed Lagrangian $(k+1)$-handle which
enables us to build the cobordism $V$, as well as the idea of
implanting it along an isotropic disc, is inspired by a construction
of Dimitroglou-Rizell appearing in \cite{DimRizell--LegAmbSurg}, which
implements $k$-surgery for Legendrian submanifolds and builds
corresponding Lagrangian cobordisms (in a different sense of the
word, see Section \ref{sec:relat-constr-DR}).

\subsection{Lagrangian cobordisms}
\label{sec:lagrangian-cobordisms}
The notion of Lagrangian cobordism appearing in Theorem
\ref{main-thm:antisurgery} is that of Biran--Cornea
\cite{Biran-Cornea--Lag-Cob-I}, adapted to the immersed setting in an
obvious way:
Two ordered collections $(\iota_i: L_i \to M)_{i=1}^r$, $(\iota_j':
L_j' \to M)_{j=1}^s$ of immersed Lagrangian submanifolds of $M$ are
called \emph{Lagrangian cobordant} if there exists a smooth cobordism
$(V;\coprod_i L_i,\coprod_j L_j')$ together with a Lagrangian
immersion $V \to [0,1] \times \mathbb R \times M \subset T^*\mathbb R
\times M$ such that for some small $\delta > 0$, we have
\begin{equation*}
    V\vert_{[0,\delta) \times \mathbb R} = \coprod_{i=1}^{r}
    ~ [0,\delta) \times \{i\} \times L_i \quad \text{and} \quad V\vert_{(1-\delta,1]\times \mathbb R} = \coprod_{j=1}^{s} ~
    (1-\delta,1]\times \{j\} \times L_j'.
\end{equation*}
Here we use the notation $V \vert_U := V \cap (U \times M)$ to denote
the part of $V$ that lies over some subset $U \subset T^*\mathbb R$,
and we identify $T^*\mathbb R \cong \mathbb R \times \mathbb R$ in the
standard way. The Lagrangian submanifold
$V \subset T^*\mathbb R \times M$ is called an \emph{immersed
    Lagrangian cobordism} with negative ends $(L_i)_{i=1}^r$ and
positive ends $(L_j)_{j=1}^r$, and this relationship is denoted by
$V: (L_1',\dots,L_s') \leadsto (L_1,\dots,L_r)$. In this article we
will only deal with the case $r = s = 1$, i.e., with Lagrangian
cobordisms
\begin{equation*}
    V: L' \leadsto L
\end{equation*}
with a single positive and a single negative end.

Lagrangian cobordisms have recently attracted a lot of interest due to
the fact that, provided certain monotonicity assumptions hold, they
preserve Floer theoretic invariants and encode information about the
Fukaya category, see
\cite{Biran-Cornea--Lag-Cob-I,Biran-Cornea--Lag-Cob-II} and also the
recent \cite{Mak-Wu--Dehn-twists}. So far, there have been essentially
two known constructions of Lagrangian cobordisms, which are based on
Hamiltonian isotopy resp.\ Lagrangian 0-surgery. Extending the toolkit
for building new ones, such as those provided by Theorem
\ref{main-thm:antisurgery} and Theorem \ref{main-thm:embedded-cob}
below, was one of the motivations for the present
paper.

\subsection{Desingularization}
\label{sec:surgering-again}
The newly created double point $p$ of the Lagrangian $L'$ resulting
from Lagrangian antisurgery can be resolved by Lagrangian 0-surgery.
Provided that $L$ is embedded, this yields a Lagrangian $L^\natural$
which is also embedded and diffeomorphic to $L' \# P^n$ or to
$L' \# Q^n$, where $P^n = S^1 \times S^{n-1}$, and $Q^n$ is the
mapping torus of an orientation-reversing involution of $S^{n-1}$.

There are in fact two families of such resolutions which correspond to
the two ways ordering the sheets meeting at $p \in L'$ (as is always
the case for Lagrangian 0-surgery). We will show that there exists one
such family such that for all $L^\natural$ in that family which are of
sufficiently small size (in the sense of Definition
\ref{defn:size-of-0-surgery}), one can in fact extend the resolution
such as to simultaneously remove the singular locus of the immersed
Lagrangian cobordism $V: L' \leadsto L$ produced by Theorem
\ref{main-thm:antisurgery}:

\begin{thm}
    \label{main-thm:embedded-cob}
    There exists a choice of ordering of the sheets meeting at the
    double point $p \in L'$ such that for all $L^\natural$ in the
    corresponding family of resolutions of $p$ which are of
    sufficiently small size, there exists an \emph{embedded}
    Lagrangian cobordism $V^\natural: L^\natural \leadsto L$ which
    coincides with the immersed cobordism $V: L' \leadsto L$ outside
    of a small neighbourhood of the singular locus of $V$. As a smooth
    manifold, $V^\natural$ is diffeomorphic to the manifold obtained
    from $[0,1] \times L$ by consecutively attaching a $(k+1)$-handle
    and a 1-handle.
\end{thm}

The construction which constitutes the proof of Theorem
\ref{main-thm:embedded-cob} will be given in Section
\ref{sec:singular-locus}. As stated in Theorem
\eqref{main-thm:antisurgery}, the singular locus of $V$ looks like a
line of double points (see Proposition \ref{prop:sing-locus} for a
precise description of the singularity of the correponding model
the passage from $V$ to $V^\natural$ replaces a neighbourhood of it by
a Lagrangian 1-handle.

One should note at this point that our construction of antisurgery
cobordisms cannot be replaced by simply appealing to the $h$-principle
\cite{eliashberg-mishachev02} satisfied by immersed Lagrangian
cobordisms: Indeed, this would not provide the amount of information
about the singular locus that we need in order to control the topology
of the result of desingularizing the immersed cobordism by a version
of Lagrangian surgery.

\subsection{Reversing the construction}
\label{sec:revers-constr}
Lagrangian antisurgery constructs from a Lagrangian $L \subset M$ an
new Lagrangian $L'$ with one (additional) double point. Changing
perspectives, we can view $L$ as the result of resolving a double
point of $L'$ by an operation which is an $(n-k-1)$-surgery on the
level of smooth manifolds, and which we therefore refer to as
\emph{Lagrangian $(n-k-1)$-surgery}. We will discuss this
point of view in Section \ref{sec:reverse-constr}.

In the case $k = n-1$, this reversed operation is the same as
classical Lagrangian $0$-surgery, up to Lagrangian isotopy. That is,
if $L'$ is the result of an $(n-1)$-antisurgery on $L$, then $L$ can
be obtained back from $L'$ by classical Lagrangian $0$-surgery
followed by a Lagrangian isotopy, and vice versa. As a consequence,
for $k=n-1$ the ends of the desingularized antisurgery cobordisms
$V^\natural: L^\natural \leadsto L$ are both resolutions of the
Lagrangian $L'$ by Lagrangian 0-surgery. Recall that Theorem
\ref{main-thm:embedded-cob} only asserts the existence of such a
cobordism for $L^\natural$ belonging to \emph{one} of the two families
of resolutions $p \in L'$. The next proposition identifies which one
it is:

\begin{prop}[See Proposition \ref{prop:distinct_0_surgeries}]
    \label{main-prop:distinct_0_surgeries}
    The ends of the desingularized cobordism
    $V^\natural: L^\natural \leadsto L$ resulting from
    $(n-1)$-antisurgery on $L$ belong to \emph{distinct} families of
    resolutions of $p \in L'$ by Lagrangian 0-surgery.
\end{prop}

\subsection{Cobordisms between Clifford and Chekanov tori}
\label{sec:monot-outp-surg}
As an application of Theorem \ref{main-thm:embedded-cob}, we will
construct cobordisms between Clifford tori $T^2_{Cl}(a)$ and Chekanov
tori $T^2_{Ch}(A)$ in $\mathbb R^4$. These are monotone Lagrangian
tori which, in both cases, are specified uniquely up to Hamiltonian
isotopy by the areas $a, A > 0$ of any disc of Maslov index 2 with
boundary on them.

\begin{thm}[See Theorem \ref{thm:Clifford-Chekanov}]
    \label{main-thm:Clifford-Chekanov}
    For every choice of $a < A$ with $a/A$ sufficiently close to 1,
    there exists a Lagrangian cobordism
    $T^2_{Cl}(a) \leadsto T^2_{Ch}(A)$ which as smooth a manifold is
    obtained from $[0,1] \times T^2$ by successively attaching a
    2-handle and a 1-handle.
\end{thm}

To put Theorem \ref{main-thm:Clifford-Chekanov} into context, recall
first the classical fact that there does not exist a Hamiltonian
isotopy between any two Chekanov and Clifford tori
\cite{Chekanov--Lagr-cobs, eliashberg-polterovich1997}. Since the
relation of being Lagrangian cobordant is a generalization of the
relation of being Hamiltonian isotopic (as Hamiltonian isotopies give
rise to \emph{Lagrangian suspension} cobordisms $L \leadsto \phi_t(L)$, see
\cite{Biran-Cornea--Lag-Cob-I}), this fact can also be viewed as a
restriction on the type of Lagrangian cobordisms that can exist
between Clifford and Chekanov tori.

One way of disproving the existence of a Hamiltonian isotopy between
$T^2_{Cl}(A)$ and $T^2_{Ch}(A)$ is to note that the numbers of
pseudoholomorphic discs of Maslov index 2 through a generic points on
these tori are different (this number is 1 for the Chekanov torus, but
2 for the Clifford torus), while the existence of a Hamiltonian
isotopy between them would imply that these counts are the same. A
similar argument precludes the existence of any \emph{monotone}
Lagrangian cobordism between $T^2_{Cl}(A)$ and $T^2_{Ch}(A)$, because
such cobordisms also preserve counts of Maslov 2 discs
\cite{Chekanov--Lagr-cobs,Biran-Cornea--Lag-Cob-I}. In particular, the
statement of Theorem \ref{main-thm:Clifford-Chekanov} cannot be
extended to include the case $a = A$, because the resulting cobordism
would be automatically monotone (see Proposition
\ref{prop:H_2-surjectivity}).

In fact, one can adapt the argument used to prove the latter to obtain
the following statement: There does not exist a Lagrangian cobordism
$V: T^2_{Cl}(A) \leadsto T^2_{Ch}(a)$ with $a < A$ (i.e., connecting a
Clifford torus to a \emph{smaller} Chekanov torus, in contrast to
Theorem \ref{main-thm:Clifford-Chekanov}) and with the property that
\begin{equation}
    \label{eq:min-area-prop}
    \inf \,\big\{\omega(\sigma)~|~ \sigma \in \pi_2(T^*\mathbb R \times
    \mathbb R^4, V), ~\omega(\sigma) > 0 \big\} = a.
\end{equation}
Indeed, if one assumes that such a cobordism exists, one arrives at a
contradiction when considering the number of Maslov 2 discs with
boundary of the cobordism passing through a generic point and
representing the push-forward of the unique class in
$H_2(\mathbb R^4, T_{Ch}^2(a))$ represented by a Maslov 2
pseudoholomorphic disc: On the one hand, this number would have to be
independent of the chosen point, as property \eqref{eq:min-area-prop}
ensures that there can be no bubbling and hence that the corresponding
moduli space is compact. On the other hand, it would need to be 1 for
a point on the Chekanov end, but 0 for a point on the Clifford end.

In contrast to that, there is no obvious obstruction to the existence
of a cobordism $T^2_{Cl}(a) \leadsto T^2_{Ch}(A)$ with property
\eqref{eq:min-area-prop}, i.e., from a Clifford torus to a
\emph{larger} Chekanov torus, as compactness of the moduli space of
discs described is not guaranteed in this situation.

To connect these observations with Theorem
\ref{main-thm:Clifford-Chekanov}, note that a cobordism between
Lagrangian 2-tori whose topology is that of a 1-antisurgery cobordism,
i.e., as described in the theorem, has property
\eqref{eq:min-area-prop} if $A = k \cdot a$ for some $k \in \mathbb N$
(this follows e.g.\ from Proposition
\ref{prop:H_2-surjectivity}). Therefore, it seems plausible to
conjecture that the r\^ole of Clifford and Chekanov in Theorem
\ref{main-thm:Clifford-Chekanov} cannot be swapped, i.e., that a
cobordism of this topology connecting a Clifford torus $T_{Cl}^2(A)$
to a \emph{smaller} Chekanov torus $T_{Ch}^2(a)$ does not exist. (Note
however that this is just a conjecture, as the existence of such a
cobordism for $a/A$ close to 1 is not ruled out by the arguments in
the previous paragraph). On the other hand, it is likely that one can
extend the statement of Theorem \ref{main-thm:Clifford-Chekanov} to
give the existence of an antisurgery cobordism
$T_{Cl}^2(a) \leadsto T_{Ch}^2(A)$ for \emph{any} choice of
$0 < a < A$ by deforming the local model for antisurgery suitably
(cf.\ Section \ref{sec:cobordisms-clifford-chekanov}).

While the quest for the ``simplest'' Lagrangian cobordism connecting
given Clifford and Chekanov tori is interesting and subtle, we note
that the existence problem for such tori is completely flexible if one
does not constrain the topology of the cobordisms one considers: As
mentioned above, \emph{immersed} Lagrangian cobordisms are governed by
an $h$-principle, and an immersed cobordism can always be turned into
an embedded one by removing all transverse double points by Lagrangian
0-surgery (which are the only singularities after applying a small
perturbation).

\subsection{Relation to other work}
\label{sec:relat-prev-work}
As mentioned before, one important source of inspiration for our
construction is the surgery construction for Legendrian submanifolds
appearing in \cite{DimRizell--LegAmbSurg}. The local model for the
immersed Lagrangian handle we use can be traced back to
\cite{Arnold--Lagr-Leg-Cob-I, Audin-Lalonde-Polterovich}, where it
appears in a slightly different guise. It seems that the passage from
$L$ to $L^\natural$ in the case $n=2$ and $k=1$ is identical to an
operation described in \cite{Yau--SurgeryLagrSurfaces2013}. The
article \cite{Mak-Wu--Dehn-twists} is an exploration of relations
between Lagrangian surgery and Lagrangian cobordisms in a different
direction. Cobordisms as described in Theorem
\ref{main-thm:Clifford-Chekanov} have recently also been constructed
by Jeff Hicks \cite{Hicks--Lagr-cobs-and-wallcrossing} using Lefschetz
fibrations.

\section{Lagrangian antisurgery}
\label{sec:constr-cobord}
In this section we will explain the construction of immersed
Lagrangian $(k+1)$-handles $\Gamma \subset T^*\mathbb R \times
T^*\mathbb R^n$ for $0 \leq k \leq n-1$, which will serve as the
local models for the construction of the cobordisms appearing in
Theorem \ref{main-thm:antisurgery}. Theses handles are immersed
Lagrangian cobordisms
\begin{equation*}
    \Gamma: \Lambda' \leadsto \Lambda
\end{equation*}
diffeomorphic to $D^{k+1} \times D^{n-k}$ and whose ends are
Lagrangian submanifolds $\Lambda \approx S^k \times D^{n-k}$ and
$\Lambda' \approx D^{k+1} \times S^{n-k-1}$ of $T^*\mathbb R^n$. The
construction is inspired by a similar one in
\cite{DimRizell--LegAmbSurg} (see Section \ref{sec:relat-constr-DR}
for the precise relationship).

Throughout, we will use the standard symplectic form on
$T^*\mathbb R \times T^*\mathbb R^n$ given by
\begin{equation*}
    \omega_{\mathrm{std}} = dx_0 \wedge dy_0 + \sum_{i=1}^n dx_i \wedge dy_i,
\end{equation*}
where $x_0, y_0$ and $x_1, \dots, x_n, y_1, \dots, y_n$ are the usual
coordinates on $T^*\mathbb R$ resp.\ $T^* \mathbb R^n$.

\subsection{Construction of $\Gamma$}
\label{sec:immersed-cobordism} 
The handle $\Gamma$ will be defined as the union of the graphs of
exact 1-forms $+dF$ and $-dF$, where $F: \mathfrak U \to \mathbb R$ is a
function defined on a certain subset $\mathfrak U \subset \mathbb R \times
\mathbb R^n$.

As a first step in defining $\mathfrak U$ and $F$, consider smooth
functions $\sigma: \mathbb R \to \mathbb R_{\geq 0}$ and $\rho:
\mathbb R_{\geq 0} \to \mathbb R_{\geq 0}$ satisfying
\begin{enumerate}
\item $\sigma(x_0) = 0$ for $x_0 \leq \delta$,
\item $\sigma(x_0) = 1+\varepsilon$ for $x_0 \geq 1 - \delta $,
\item $\sigma'(x_0) > 0$ for $\delta < x_0 < 1 - \delta$,
\end{enumerate}
and
\begin{enumerate}
\item $\rho(r^2) = 1$ for $r^2$ close to 0,
\item $\rho(r^2) = 0$ for $r^2 \geq 1 + 2 \varepsilon$,
\item $-1/(1+\varepsilon) < \rho'(r^2) \leq 0$ for all $r^2 \in
    \mathbb R_{\geq 0}$
\end{enumerate}
for certain small constants $\varepsilon, \delta > 0$. Denote by $r^2,
s^2: \mathbb R^n \to \mathbb R_{\geq 0}$ the functions given by
$r^2(\x) = x_1^2 + \dots + x_{k+1}^2$ and $s^2 (\x) = x_{k+2}^2 +
\dots + x_n^2$, where $\x = (x_1,\dots,x_n) \in \mathbb R^n$. Then
define a function $f: \mathbb R \times \mathbb R^n \to \mathbb R$ by
\begin{equation}
    f(x_0,\x) = r^2 + \sigma(x_0)\rho(r^2) -s^2 -1
\end{equation}
for $(x_0,\x) \in \mathbb R \times \mathbb R^n$, where $r^2 \equiv
r^2(\x)$ and $s^2 \equiv s^2(\x)$.

\begin{figure}[t]
    \centering
    \includegraphics[scale=0.85]{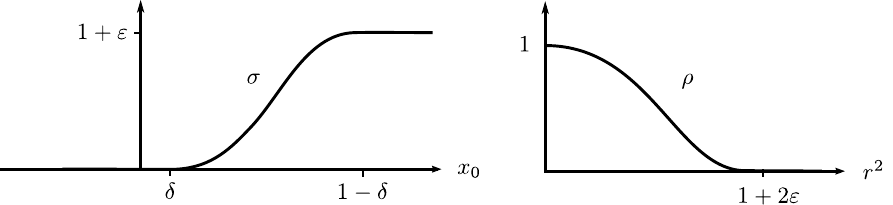}
    \caption{The auxiliary functions
        $\sigma: \mathbb R \to \mathbb R_{\geq 0}$ and
        $\rho: \mathbb R_{\geq 0} \to \mathbb R_{\geq 0}$ used in the
        construction of the Lagrangian handles $\Gamma$ which will
        serve as the local models for the Lagrangian antisurgery
        cobordisms.}
\end{figure}

Consider now the set $\mathfrak U = \{(x_0,\x) \in \mathbb R \times \mathbb R^n ~|~
f(x_0,\x) \geq 0\}$ and define $F: \mathfrak U \to \mathbb R$ by
\begin{equation}
    \label{eq:1}
    F(x_0,\x) = f(x_0,\x)^{3/2}.
\end{equation}
The restriction of $F$ to $\mathrm{int}(\mathfrak U)$ is smooth, with
differential given by
\begin{equation}
    \label{eq:dF}
    dF = \frac{3}{2}f(x_0,\x)^{1/2}\big(\sigma'(x_0)\rho(r^2)dx_0 + (1 + \sigma(x_0)\rho'(r^2))dr^2 -ds^2\big).
\end{equation}
Note that $dF$ extends to a section of
$T^*\mathbb R \times T^*\mathbb R^n$ defined on all of $\mathfrak U$
which vanishes along $\partial \mathfrak U = \{f(x_0,\x) = 0\}$; we
will denote this extended section by $dF$ as well. The graphs
$\Gamma_\pm \subset T^*\mathbb R \times T^*\mathbb R^n$ of
$\pm dF: \mathfrak U \to T^*\mathbb R \times T^*\mathbb R^n$ are
Lagrangian submanifolds with boundary, and the tangent spaces along
the boundary are given by
\begin{equation*}
    T \Gamma_\pm \vert_{\partial \Gamma_\pm} = T( N^*(\partial \mathfrak U))\vert_{\partial \mathfrak U},
\end{equation*}
where $N^*(\partial \mathfrak U)$ denotes the conormal bundle of $\partial
\mathfrak U$
. Hence $\Gamma_+$ and $\Gamma_-$ fit together smoothly along
$\partial \mathfrak U$, in the sense that their union
\begin{equation}
    \label{eq:Gamma}
    \Gamma = \Gamma_+ \cup \Gamma_- = \{((x_0,\x),\pm dF(x_0,\x))~|~(x_0,\x) \in \mathfrak U\}
\end{equation}
is a submanifold of $T^*\mathbb R \times T^* \mathbb R^n$ which is
embedded near $\partial \mathfrak U$. The singular locus
$\Gamma^s \subset \Gamma$ along which $\Gamma$ is not embedded is
the set of points $(x_0,\x) \in \mathrm{int}\,\mathfrak U$ at which
$dF$ vanishes, which is given by
\begin{equation}
    \label{eq:2}
    \Gamma^s = \{((x_0, 0), (0, 0)) \in T^*\mathbb R \times
    T^*\mathbb R^n ~|~ x_0 \geq 1 - \delta\},
\end{equation}
see Proposition \ref{prop:sing-locus}; that is, $\Gamma^s$ is a
1-dimensional family of double points.

$\Gamma$ is the immersed image of a $(k+1)$-handle $D^{k+1} \times
D^{n-k}$
, and moreover an immersed Lagrangian cobordism
\begin{equation*}
    \Gamma: \Lambda'\leadsto \Lambda
\end{equation*}
in the sense of Section \ref{sec:lagrangian-cobordisms}. To see the
latter and to describe the ends, set
$\mathfrak U_{x_0} = \{\x \in \mathbb R^n ~|~ (x_0,\x) \in \mathfrak
U\}$ for $x_0 \in \mathbb R$ and define
$F_{x_0}: \mathfrak U_{x_0} \to \mathbb R$ to be the function given by
$F_{x_0}(\x) = F(x_0,\x)$ for $\x \in \mathfrak U_{x_0}$.  Since
$F_{x_0}$ is independent of $x_0$ if either $x_0 \leq \delta$ or
$x_0 \geq 1 -\delta$, it follows that the part of $\Gamma$ lying over
$ (-\infty,\delta] \times \mathbb R \cup [1-\delta,\infty) \times
\mathbb R \subset T^*\mathbb R$ is
\begin{equation}
    \label{eq:17}
    (-\infty,\delta] \times \{0\} \times  \Lambda \quad \cup \quad [1
    -\delta,\infty) \times \{0\} \times \Lambda'
\end{equation}
with
\begin{equation}
    \label{eq:Ends}
    \begin{aligned}
        \Lambda &= \{(\x,\pm dF_0(\x)) \in T^*\mathbb R^n~|~ \x \in
        \mathfrak U_0\}, \\
        \Lambda' &= \{(\x,\pm dF_1(\x)) \in T^*\mathbb
        R^n~|~ \x \in \mathfrak U_1\}.
    \end{aligned}
\end{equation}
This shows that $\Gamma$ is a Lagrangian cobordism (up to modifying
the ends in an obvious way).

\subsection{Isotropic surgery discs}
\label{sec:surgery-disc}
We will now describe the situation in which it is possible to implant
the local model described above, such as to produce from a given
Lagrangian $L$ a new immersed Lagrangian $L'$ together with a
Lagrangian cobordism $V: L' \leadsto L$. 

The following definition is an adaptation of Definition 4.2 in
\cite{DimRizell--LegAmbSurg} to our setting.

\begin{defn}
    \label{defn:surgery-disc}
    Let $L \subset M$ be a Lagrangian submanifold and let $S \subset L$ be
    a embedded $k$-sphere with trivializable normal bundle. An
    \emph{isotropic surgery disc} for $S$ consists of the following data:
    \begin{enumerate}
    \item An embedded isotropic $(k+1)$-disc $D \subset M$ such
        that
        \begin{itemize}
        \item $\partial D = S$,
        \item $\mathrm{int}\,D \cap L = \emptyset$,
        \item any vector field $X \subset TD\vert_S$ which is
            outward-pointing normal to $S = \partial D$ is nowhere
            contained in $TL$.
        \end{itemize}
    \item A symplectic subbundle $E$ of $(TD)^\omega$ such that $TD
        \oplus E = (TD)^\omega$
        , and a symplectic trivialization $\Psi: D \times \mathbb
        C^{n-k-1} \to E$ such that the Lagrangian subbundle $\Psi(S
        \times \mathbb R^{n-k-1})$ of $E \vert_S$ is contained in $TL
        \vert_S$.
    \end{enumerate}
    We will usually denote isotropic surgery discs simply by $D$,
    omitting the bundle $E$ and its trivialization $\Psi$ from the
    notation.
\end{defn}

An isotropic surgery disc $D \equiv (D, E, \Psi)$ for a sphere
$S \subset L$ determines a homotopy class of trivializations of the
normal bundle of $S \subset L$ as follows: Let $Y \subset TL \vert_S$
be a vector field which is normal to $S \subset L$ and such that
$\omega(X,Y) > 0$ for a vector field $X \subset TD\vert_S$ which is
outward-pointing normal to $S$ (such a vector field $Y$ exists due to
the assumption on such $X$ in Definition
\ref{defn:surgery-disc}). Then the subbundle
$\Psi(S \times \mathbb R^{n-k-1}) \oplus \mathbb R Y$ of $TL\vert_S$
is complementary to $TS$ and of rank $n-k$, and thus it spans the
normal bundle of $S \subset L$. Since the space of all such vector
fields $Y$ is non-empty and contractible, the corresponding
trivialization is determined up to
homotopy. 

\begin{example}
    \label{ex:surgery-disc}
    The prototypical example for the situation described in Definition
    \ref{defn:surgery-disc} is given by the Lagrangian $\Lambda \subset
    T^*\mathbb R^n$ described in \eqref{eq:Ends} and the $k$-sphere
    \begin{equation}
        \label{eq:S_m}
        S_0 =
        \{(\x, \y) \in T^*\mathbb R^n ~|~ x_1^2 + \dots + x_{k+1}^2 = 1,
        x_{k+2} = \dots = x_n = 0, \y = 0\};
    \end{equation}
    the obvious choice of isotropic surgery disc for $S_0 \subset
    \Lambda$ is
    \begin{equation}
        \label{eq:D_0}
        D_0 = \{(\x,\y) \in T^*\mathbb R^n ~|~ x_1^2 + \dots + x_{k+1}^2 \leq 1,
        x_{k+2} = \dots = x_n = 0, \y = 0\}   
    \end{equation}
    together with the symplectic subbundle
    \begin{equation}
        \label{eq:18}
        E_0 = \langle \partial_{x_{k+2}}, \dots
        , \partial_{x_n},\partial_{y_{k+2}},\dots,\partial_{y_n} \rangle
    \end{equation}
    of $(TD_0)^\omega$ and the identification $\Psi_0: D_0 \times
    \mathbb C^{n-k-1} \to E_0$ taking $D_0 \times \mathbb R^{n-k-1}$
    to the subbundle
    $\langle \partial_{x_{k+2}},\dots,\partial_{x_n}\rangle$ and $D_0
    \times i\mathbb R^{n-k-1}$ to the subbundle
    $\langle \partial_{y_{k+2}},\dots,\partial_{y_n}\rangle$.

\end{example}

Assume that we are in the situation of Definition
\ref{defn:surgery-disc}, i.e., that we have a Lagrangian $L$ with a
sphere $S \subset L$ and a corresponding isotropic surgery disc $D
\equiv (D,E,\Psi)$. Let $\phi: D_0 \to D$ be a diffeomorphism;
together with the symplectic trivialization $\Psi: D \times \mathbb
C^{n-k-1} \times E$, this determines an isomorphism of symplectic
vector bundles $T^*D_0 \oplus E_0 \cong T^*D \oplus E$ (here we use
the notation of Example \ref{ex:surgery-disc}). An application of the
isotropic neighbourhood theorem then yields an extension of $\phi$ to
a symplectomorphism
\begin{equation*}
    \phi: \mathcal W_0 \to
    \mathcal W
\end{equation*}
between appropriate Darboux-Weinstein neighbourhoods $\mathcal W_0
\supset D_0$ and $\mathcal W \supset D$ of the discs in $T^*\mathbb
R^n$ resp.\ $M$, and we may assume that this extension satisfies
\begin{equation}
    \label{eq:implant}
    \phi(\Lambda \cap \mathcal W_0) = L \cap \mathcal W.
\end{equation}
To see this, note that the condition that the outward normal vector
field to $S \subset D$ is nowhere tangent to $L$ guarantees that one
can arrange $D\phi (T\Lambda \vert_{S_0}) = TL \vert_S$; after
adjusting $\phi$ by a Hamiltonian isotopy and possibly shrinking the
Weinstein neighbourhoods, one obtains \eqref{eq:implant}.

\subsection{Implantation of the local model and proof of Theorem \ref{main-thm:antisurgery}}
\label{sec:defn-antisurgery}
We now explain how to implant the local model and give the definition
of Lagrangian antisurgery and of the corresponding Lagrangian
cobordism. This will complete the proof of Theorem
\ref{main-thm:antisurgery} (up to the description of the singular
locus of the cobordism, for which we refer to Proposition
\ref{prop:sing-locus}).

As before, we assume that we have a Lagrangian $L$ containing an
embedded sphere $S \subset L$ together with an isotropic surgery disc
$D$. To prepare the construction, consider the neighbourhood of
$D_0 \subset T^*\mathbb R^n$ (see Example \ref{ex:surgery-disc}) given
by
\begin{equation}
    \label{eq:U_m}
    \mathcal U_0 = \Big\{(\x,\y) \in T^*\mathbb R^n ~|~
    r^2 < 1 + 2 \varepsilon,~  s^2 < 2 \varepsilon, ~
    \Vert \y \Vert^2 < 6\sqrt{2 \varepsilon}(1 + 4 \varepsilon)\Big\}
\end{equation}
and denote by $\mathcal U_0^c$ the complement of $\mathcal U_0$ in
$T^*\mathbb R^n$.  
The following lemma shows that the part of the model cobordism
$\Gamma$ that projects to $\mathcal U_0^c \subset T^*\mathbb R^n$ 
lies over $\mathbb R \subset T^*\mathbb R$ and is ``cylindrical'':

\begin{lem}
    \label{lem:cyl-part}
    We have $\Gamma \cap (T^*\mathbb R \times \mathcal U_0^c) =
    \mathbb R \times (\Lambda \cap \mathcal U_0^c) = \mathbb R \times
    (\Lambda' \cap \mathcal U_0^c)$. 
\end{lem}
\begin{proof}
    We first claim that for
    $((x_0,\x),(y_0,\y)) = ((x_0,\x),\pm dF(x_0,\x)) \in \Gamma$ with
    $r^2 < 1 + 2 \varepsilon$, we already have
    $(\x,\y) \in \mathcal U_0$. To see that, recall that the set
    $\mathfrak U \subset \mathbb R \times \mathbb R^n$ over which
    $\Gamma$ lives is characterized by $f(x_0,\x) \geq 0$, where
    $f(x_0,\x) = r^2 + \sigma(x_0)\rho(r^2) - s^2 - 1$. Since
    $r^2 \mapsto r^2 + \sigma(x_0)\rho(r^2) - 1$ is strictly
    increasing (this follows from the assumptions that
    $\sigma(x_0) \leq 1 + \varepsilon$ for all $x_0$ and
    $-1/(1 + \varepsilon) < \rho'(r^2)$ for all $r^2$) with value
    $2 \varepsilon$ at $r^2 = 1 + 2 \varepsilon$ for every
    $x_0 \in \mathbb R$, it follows that $s^2 < 2
    \varepsilon$. Moreover, one can read off from the expression
    \eqref{eq:dF} for $dF(x_0,\x)$ that the bound on
    $\Vert \y \Vert^2$ is satisfied whenever $r^2 < 1 + 2 \varepsilon$
    and $s^2 < 2 \varepsilon$.
    Let now
    $((x_0,\x),(y_0,\y)) \in \Gamma \cap (T^*\mathbb R \times \mathcal
    U_0^c)$. As a consequence of the claim above, we obtain
    $r^2 \geq 1 + 2 \varepsilon$, and hence \eqref{eq:dF} for
    $dF(x_0,\x)$ simplifies to
    $dF(x_0,\x) = \frac{3}{2} (r^2 - s^2 -1)^{1/2}(dr^2 - ds^2)$, as
    $\rho(r^2) \equiv 0$ for $r^2 \geq 1 + 2 \varepsilon$. Since this
    has vanishing $dx_0$ component $y_0 = 0$ and is independent of
    $x_0$, it follows that
    $((x_0,\x),(y_0,\y)) = ((x_0, \x), (0, \y))$ lies in
    $\mathbb R \times (\Lambda \cap \mathcal U_0^c)$ and in
    $\mathbb R \times (\Lambda' \cap \mathcal U_0^c)$. Thus
    $\Gamma \cap (T^*\mathbb R \times \mathcal U_0^c)$ is contained in
    both of these sets. The inclusions in the other direction are
    obvious.
\end{proof}

The neighbourhood $\mathcal U_0$ of $D_0$ can be made arbitrarily
small by letting the parameter $\varepsilon$ tend to zero. In
particular, by choosing the parameter $\varepsilon$ sufficiently
small, we may assume that the closure $\overline{\mathcal U_0}$ is
contained in a Weinstein neighbourhood $\mathcal W_0$ of $D_0$ as
described in Section \ref{sec:surgery-disc}, i.e., such that we have a
symplectic identification $\phi: \mathcal W_0 \to \mathcal W$ with a
Weinstein neighbourhood $\mathcal W$ of $D$.

\begin{defn}
    \label{defn:surgery}
    Given such choices of $\varepsilon$ and $\phi$, we define the
    immersed Lagrangian $L' \subset M$ obtained from $L$ by
    \emph{Lagrangian $k$-antisurgery} along the isotropic disc $D$ by
    \begin{equation}
        \label{eq:defn-antisurgery}
        L' = (L \cap \mathcal W^c) ~ \cup ~ \phi(\Lambda' \cap \mathcal W_0),
    \end{equation}
    and its \emph{immersed Lagrangian trace} $V: L' \leadsto L$ by
    \begin{equation}
        \label{eq:defn-cob}
        V = \mathbb R \times (L \cap \mathcal W^c) ~
        \cup ~(\mathrm{id} \times \phi)(\Gamma \cap (T^*\mathbb R
        \times \mathcal W_0)),
    \end{equation}
    using the symplectomorphism $\mathrm{id} \times \phi: T^*\mathbb R
    \times \mathcal W_0 \to T^*\mathbb R \times \mathcal W$.
\end{defn}

The fact that the pieces which we glue in fit together as required is
a consequence of Lemma \ref{lem:cyl-part}, which implies that
$\phi(\Lambda' \cap (\mathcal W_0 \smallsetminus \overline{\mathcal
    U_0})) = L \cap (\mathcal W \smallsetminus \overline{\mathcal U})$
and
$(\mathrm{id} \times \phi)(\Gamma \cap (T^*\mathbb R \times (\mathcal
W_0 \smallsetminus \overline{\mathcal U_0}))) = \mathbb R \times (L
\cap (\mathcal W \smallsetminus \overline{\mathcal U}))$, where
$\mathcal U = \phi(\mathcal U_0)$; hence the pieces of $\Lambda'$
resp.\ $\Gamma$ we glue overlap with corresponding pieces of $L$
resp.\ $\mathbb R \times L$ as required.

The Lagrangian submanifold $L' \subset M$ given by Definition
\ref{defn:surgery} is the immersed image of the manifold obtained from
$L$ by a $k$-surgery along $S$ with respect to the trivialization of
the normal bundle of $S \subset L$ determined by the surgery disc
$D$. The Lagrangian cobordism $V: L' \leadsto L$ is the immersed image
of the trace corresponding to that surgery.

\subsection{Relation to the construction in \cite{DimRizell--LegAmbSurg}}
\label{sec:relat-constr-DR}
The Lagrangian handle $\Gamma$ constructed in this section is closely
related to the Lagrangian handle constructed by Dimitroglou-Rizell in
\cite[Section 4]{DimRizell--LegAmbSurg}. To make the connection
precise, consider the function
\begin{equation*}
    \hat F: \mathfrak U \to \mathbb R, \quad (x_0,\x) \mapsto
    F(x_0,\x)\cdot x_0,
\end{equation*}
where $\mathfrak U \subset \mathbb R \times \mathbb R^n$ and
$F: \mathfrak U \to \mathbb R$ are as defined in Section
\ref{sec:immersed-cobordism}. The handle in
\cite{DimRizell--LegAmbSurg} is obtained by gluing the graphs
$\hat \Gamma_{\pm} \subset T^*(\mathbb R \times \mathbb R^n)$ of
$\pm d\hat F: \mathfrak U \to T^*(\mathbb R \times \mathbb R^n)$, in
analogy with the construction of $\Gamma$ in \ref{eq:Gamma}, such as
to obtain
\begin{equation*}
    \hat \Gamma = \hat \Gamma_- \cup \hat \Gamma_+ = \{((x_0,\x),\pm
    d\hat F(x_0,\x))~|~(x_0,\x) \in \mathfrak U\}.
\end{equation*}

The part of $\hat \Gamma$ lying over $\mathfrak U \cap \{x_0 > 0 \}$
is then mapped to $\mathbb R \times J^1(\mathbb R^n)$, the
symplectization of the 1-jet space of $\mathbb R^n$, using a
symplectic identification
$T^*(\mathbb R_{> 0} \times \mathbb R^n) \cong \mathbb R \times
J^1(\mathbb R^n)$, which results in a Lagrangian cobordism with ends
which are cylindrical over \emph{Legendrian} submanifolds of
$J^1(\mathbb R^n)$ (see e.g.\ the introduction of
\cite{DimRizell--LegAmbSurg} for the relevant definitions). In
particular, the result is not a Lagrangian cobordism in the sense of
\cite{Biran-Cornea--Lag-Cob-I} (and neither is $\hat \Gamma$), as it
does not have the required cylindrical ends described in Section
\ref{sec:lagrangian-cobordisms}.

\section{Lagrangian antisurgery and surgery}
\label{sec:reverse-constr}
In this section we explain the relationship between Lagrangian
$(n-1)$-antisurgery and classical Lagrangian $0$-surgery
\cite{Polterovich--Surgery-of-Lagrangians}. We start by recalling the
construction of Lagrangian 0-surgery as described e.g.\ in
\cite[Section 6.1]{Biran-Cornea--Lag-Cob-I}.

\subsection{Classical Lagrangian 0-surgery}
\label{sec:lagrangian-0-surgery}
\begin{figure}[t]
    \centering
    \includegraphics[scale=.95]{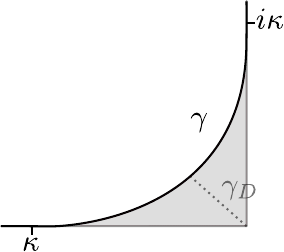}
    \caption{A curve $\gamma$ of the type used in the definition of Lagrangian 0-surgery.}
    \label{fig:surgery-curve}
\end{figure}

Let $\gamma = (a,b): \mathbb R \to T^*\mathbb R$ be an embedded smooth
curve satisfying
\begin{align}
    \label{eq:surgery-curve-spec}
  \begin{split}
      \gamma(t) &= (t,0),~ t \in (-\infty,-\kappa],\\
      a(t) &< 0 < b(t),~ t \in (-\kappa,\kappa),\\
      \gamma(t) &= (0,t),~ t \in [\kappa,\infty),
  \end{split}
\end{align}
where $\kappa > 0$ is a small parameter, see Figure
\ref{fig:surgery-curve}. Then consider the embedding
\begin{equation}
    \label{eq:surgery-profile}
    h_\gamma: \mathbb R \times S^{n-1} \to T^*\mathbb R^n, \quad (t,\x)
    \mapsto (a(t) \x, b(t) \x)
\end{equation}
where $S^{n-1} = \{\x \in \mathbb R^n ~|~ \Vert \x \Vert_2 = 1\}$ is
the unit sphere in $\mathbb R^n$. The image of $h_\gamma$ is an
embedded Lagrangian submanifold of $T^*\mathbb R^n$ which outside of
the ball $B^{2n}_{\kappa}$ of radius $\kappa$ centered at
$0 \in T^*\mathbb R^n$ coincides with
$\mathbb R^n \times \{0\} \cup \{0\} \times \mathbb R^n$.  We can also
view this Lagrangian as the orbit of $\gamma$, viewed now as living in
$ T^*\mathbb R \times \{0\} \subset T^*\mathbb R^n$, under the
$SO(n)$-action on $T^*\mathbb R^n$ given by
\begin{equation}
    \label{eq:SO(n)-action}
    A(\x,\y) = (A\x, A\y)
\end{equation}
for $A \in SO(n)$.

Given now a Lagrangian $L \subset M$ with a transverse double point
$p \in L$, we can implant this local model using a Darboux chart which
identifies neighbourhoods of the sheets of $L$ meeting at $p$ with
neighbourhoods of $0$ in $\mathbb R^n \times \{0\}$ resp.\ in
$\{0\} \times \mathbb R^n$. The result is a new Lagrangian submanifold
\begin{equation}
    L^\sharp \subset M.
\end{equation}
For a fixed choice of Darboux chart, any two choices of $\gamma$
subject to the specification in \eqref{eq:surgery-curve-spec} are
related by isotopies that are constant outside of compact sets, and
these induce Lagrangian isotopies of the corresponding versions of
$L^\sharp$. Modifying this specification so that $-\gamma \cup \gamma$
lies in the first and third quadrants of $T^*\mathbb R$, instead of
the second and fourth ones as in the description above, has the same
effect as reversing the order of the sheets in the above sense and
leads results in a second family of resolutions which are Lagrangian
isotopic to one another. Resolutions which do not belong to the same
family are usually not globally Lagrangian or even smoothly isotopic,
and sometimes even distinct as smooth manifolds (e.g., orientable in
one case, but non-orientable in the other case).

In the case that $L$ is the union of two Lagrangian submanifolds $L_-$
and $L_+$ that intersect transversely at $p \in L_- \cap L_+$, we
denote by
\begin{equation*}
    L_- \# L_+
\end{equation*}
the result of any 0-surgery resulting from implanting the local
model in such a way that locally $L_-$ gets identified with
$\mathbb R^n \times \{0\}$ while $L_+$ gets identified with
$\{0\} \times \mathbb R^n$.\footnote{This is the convention used
    e.g.\ in \cite{Biran-Cornea--Lag-Cob-I}. Other papers, such as
    \cite{Polterovich--Surgery-of-Lagrangians}, use the opposite
    convention.}

Within each family of resolutions which are Lagrangian isotopic
through such ``local'' isotopies (i.e., induced by isotopies of
0-surgery models), the obstruction to being Hamiltonian isotopic is
the difference between their sizes in the sense of the following
definition:

\begin{defn}
    \label{defn:size-of-0-surgery}
    Let $L^\natural$ be the resolution a transverse double point
    $p \in L$ obtained by implanting a local model for Lagrangian
    0-surgery as described above. We define the \emph{size} of the
    resolution to be the symplectic area between the curve $\gamma$
    and the coordinate axes in the local model (the shaded region in
    Figure \ref{fig:surgery-curve}).
\end{defn}

We note that one could give a more flexible definition of Lagrangian
0-surgery including versions with non-positive size (by removing the
requirement that $\gamma$ be contained in a quadrant); however, we
only consider surgeries of positive size.

\subsection{Lagrangian 0-surgery and $(n-1)$-antisurgery}
\label{sec:lagrangian-0-surgery-1}
We now explain the connection between Lagrangian 0-surgery and
Lagrangian $(n-1)$-antisurgery, and how to see that the two operations
are inverse to one another up to Lagrangian isotopy.

\begin{figure}[t]
    \centering
    \includegraphics[scale=1]{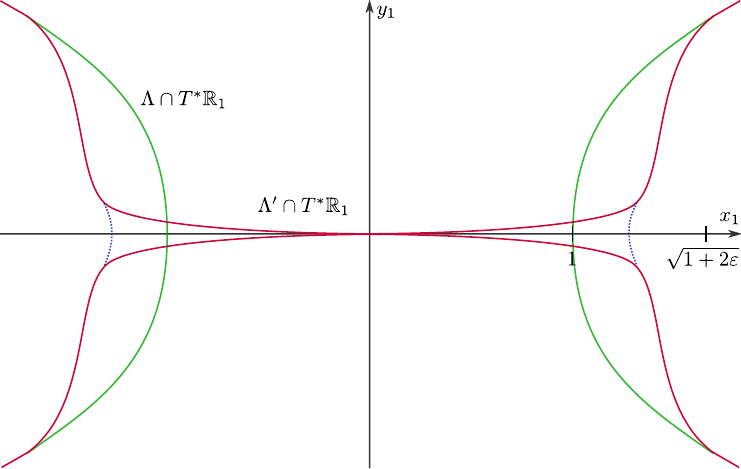}
    \caption{The figure shows the intersection of
        $\Lambda \subset T^*\mathbb R^n$, the local model for the
        non-singular end of antisurgery (in green), and
        $\Lambda' \subset T^* \mathbb R^n$, the local model for the
        singular end of antisurgery (in red), with the plane
        $T^*\mathbb R_1 \subset T^* \mathbb R^n$. By performing
        Lagrangian 0-surgery along the dotted arc, one can reverse the
        effect of $(n-1)$-antisurgery up to Lagrangian isotopy.}
    \label{fig:antisurgery-profile}
\end{figure}

Recall that Lagrangian $(n-1)$-antisurgery replaces an embedded copy
of $\Lambda \cong S^{n-1} \times D^1$ in a Lagrangian $L$ by an
immersed copy of $\Lambda' \cong D^n \times S^0$, in which the two
copies of $D^n$ intersect transversely at a double point $p$ in the
resulting Lagrangian $L'$. This double point can be resolved by
Lagrangian 0-surgery. To see that this resolution can be implemented
in such a way that the resulting Lagrangian\footnote{The
    notation $L^\flat$ is chosen to distinguish this resolution from
    $L^\natural$, the one that comes up in Theorem
    \ref{main-thm:embedded-cob} (as they belong to distinct families
    of resolutions, see Proposition \ref{prop:distinct_0_surgeries}).} $L^\flat$
is Lagrangian isotopic to $L$, we first consider the local situation
for the case $n = 1$, for which the antisurgery model is depicted in
Figure \ref{fig:antisurgery-profile}. The green arcs are the parts of
$\Lambda$ which, when performing antisurgery, are cut out and replaced
by the singular red part such as to produce $\Lambda'$. If one applies
Lagrangian 0-surgery to $\Lambda'$ in such a way that the curve
$\gamma$ in Figure \ref{fig:surgery-curve} gets mapped to the dotted
blue arc in Figure \ref{fig:antisurgery-profile}, then the resulting
non-singular Lagrangian $\Lambda^\flat$ is evidently Lagrangian
isotopic to the original $\Lambda$. Note that $\Lambda^\flat$ belongs
to the family of resolutions modelled by
\begin{equation*}
    \lambda_+ \# \lambda_-,
\end{equation*}
where $\lambda_\pm = T_{(0,0)}\Lambda'_\pm$ are the tangent spaces to
$\Lambda_\pm'= \{(\x,\pm dF_1(\x)) \in T^*\mathbb R^n~|~ \x \in
\mathfrak U_1\}$, the sheets of the singular end of the antisurgery
cobordism (cf.\ Section \ref{sec:immersed-cobordism}); this can be
read off Figure \ref{fig:antisurgery-profile}, taking into account
that $\Lambda_+'$ is the sheet which is contained in the first and
third quadrants, while $\Lambda_-$ is the sheet contained in the
second and fourth quadrants.

In order to transfer this observation to the case $n
>1$, note that the local models for both
$(n-1)$-antisurgery and 0-surgery are orbits of the respective
one-dimensional local models (viewed as living in the $(x_1,
y_1)$-coordinate subspace $T^*\mathbb R_1$ of $T^* \mathbb
R^n$), under the
$SO(n)$-action described in \eqref{eq:SO(n)-action}. One can use this
action to extend a Lagrangian isotopy between the curves $\Lambda \cap
T^* \mathbb R_1$ and $\Lambda^\flat \cap T^* \mathbb R_1$ in $T^*
\mathbb R_1$ to a Lagrangian isotopy between
$\Lambda$ and $\Lambda^\flat$ in $T^* \mathbb
R^n$ in an $SO(n)$-equivariant way.

Conversely, if we first perform Lagrangian 0-surgery at a transverse
double point $p$ of a Lagrangian $L'$ to obtain a new Lagrangian $L$,
we can reverse this operation by $(n-1)$-antisurgery after applying a
suitable Lagrangian isotopy to $L$. The Lagrangian disc required for
that can be constructed by applying the $SO(n)$-action described in
\eqref{eq:SO(n)-action} to an embedded curve
$\gamma_D: [0,1] \to T^*\mathbb R$ connecting a point on
$\mathrm{im}(\gamma)$ to $0 \in T^* \mathbb R$ in the local model for
0-surgery, see Figure \ref{fig:surgery-curve}.

\subsection{Lagrangian vs.\ Hamiltonian isotopy}
\label{sec:lagr-vs.-hamilt}
To obtain a more refined picture, and in particular to show that the
result of successively applying Lagrangian $(n-1)$-antisurgery and
0-surgery to a Lagrangian $L$ is generally not \emph{Hamiltonian}
isotopic to $L$, we will compare the symplectic areas bounded by the
ends of the local model for $(n-1)$-antisurgery.

\begin{figure}[t]
    \centering
    \includegraphics{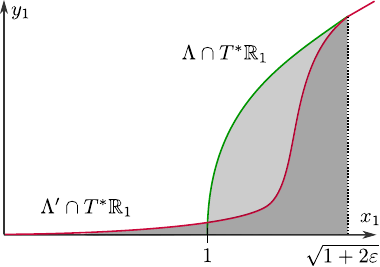}
    \caption{The areas enclosed by $\Lambda \cap T^*\mathbb R_1$
        resp.\ $\Lambda' \cap T^*\mathbb R_1$ and the line
        $\{x_1 = \sqrt{1 + 2 \varepsilon}\}$, where $\Lambda$ and
        $\Lambda'$ are the ends of model cobordism $\Gamma$
        corresponding to $(n-1)$-antisurgery.}
    \label{fig:area-difference}
\end{figure}

More precisely, we will compute the difference between the areas of
the bounded regions enclosed by the curves
$\Lambda \cap T^* \mathbb R_1$ resp.\ $\Lambda' \cap T^* \mathbb R_1$
and the line $\{x_1 = \sqrt{1 + 2 \varepsilon}\}$, i.e., twice the
difference between the two shaded areas in Figure
\ref{fig:area-difference}. The description of the ends in Section
\ref{sec:immersed-cobordism} shows that
\begin{align*}
  \Lambda \cap T^*\mathbb R_1            & = \{(x_1, \pm dF_0(x_1) \in T^*
                                 \mathbb R_1 ~|~ x_1 \geq 1\} \\
  \Lambda' \cap T^*\mathbb R_1           & = \{(x_1, \pm dF_1(x_1) \in T^* \mathbb R_1 ~|~
                                 x_1 \geq 0\},
\end{align*}
where $F_0(x_1) = (x_1^2 -1)^{3/2}$ and
$F_1(x_1) = (x_1^2 + (1 + \varepsilon)\rho(x_1^2) - 1)^{3/2}$. Since
$F_0(x_1) = F_1(x_1)$ for $x_1 \geq \sqrt{1 + 2 \varepsilon}$, the
area difference we are interested in is
\begin{align*}
  \begin{split}
      2 \left( \int_1^{\sqrt{1 + 2 \varepsilon}} \frac{\partial
              F_0}{\partial x_1} dx_1 - \int_0^{\sqrt{1 + 2
                  \varepsilon}} \frac{\partial F_1}{\partial x_1}
          dx_1\right) & = 2 \Big(F_0(\sqrt{1 + 2 \varepsilon})
          - F_0(1)\Big)\\
          & \qquad - 2 \Big(F_1(\sqrt{1 + 2 \varepsilon}) - F_1(0)\Big)
  \end{split}\\
                      & = 2 \Big((2\varepsilon)^{3/2}
                        - (2\varepsilon)^{3/2}
                        + \varepsilon^{3/2} \Big) \\
                      & = 2\varepsilon^{3/2},
\end{align*}
using that $\rho(0) = 1$ and $\rho(1 + 2 \varepsilon) = 0$. In
particular, this computation justifies that the singular end is
depicted as lying below the nonsingular end near
$x_1 = \sqrt{1 + 2 \varepsilon}$ in Figure
\ref{fig:antisurgery-profile}.

Note that an application of Lagrangian 0-surgery to resolve the double
point of the singular end $\Lambda'$ leads to a Lagrangian
$\Lambda^\flat$ which together with the line
$\{x_1 = 1 + 2 \varepsilon\}$ bounds \emph{less} area than the local
model $\Lambda'$ for the non-singular end. Consider now again a
Lagrangian $L^\flat$ resulting from applying successively Lagrangian
$(n-1)$-antisurgery and 0-surgery to a Lagrangian $L$. As a
consequence of this local picture, any Lagrangian isotopy from $L$ to
$L^\flat$ that comes from a Lagrangian isotopy between the respective
curves in the local model has non-vanishing flux and is therefore not
Hamiltonian.

\subsection{Lagrangian surgery of higher index}
Lagrangian antisurgery produces from a Lagrangian submanifold $L$ a
new Lagrangian submanifold $L'$ with an additional double
point. Switching the r\^oles of input and output, we can interpret $L$
as the result of an operation which resolves a singularity of $L'$ by
replacing an immersed copy of $D^{n-p} \times S^{p} \subset L'$ by an
embedded copy of $S^{n-p-1} \times D^{p+1} \subset L$. In the case
$p = 0$, i.e., when $L'$ is the result of Lagrangian
$(n-1)$-antisurgery on $L$, we have seen in Section
\ref{sec:lagrangian-0-surgery} that the inverse operation is classical
Lagrangian 0-surgery up to Lagrangian isotopy. In view of that, we
propose the following definition:

\begin{defn}
    \label{defn:Lagr-p-surgery}
    We say that $L$ can be obtained from $L'$ by \emph{Lagrangian
        $p$-surgery} if $L'$ can be obtained by applying Lagrangian
    $(n-p-1)$-antisurgery a Lagrangian $\widetilde L$ which is
    Lagrangian isotopic to $L$.
\end{defn}

\begin{figure}[t]
    \centering
    \includegraphics[scale=1]{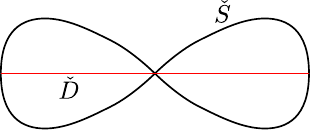}
    \caption{The central fibre of a Whitney degeneration.}
    \label{fig:Whitney-sphere}
\end{figure}

It would be interesting to characterize the admissibility of a given
immersed Lagrangian $L'$ for Lagrangian $p$-surgery by a geometric
condition analogous to the existence of an isotropic surgery disc. A
necessary condition is of course that $L'$ contains an immersed copy of
$D^{n-p} \times S^p$ which is obtained by implanting a suitable piece
of the immersed Lagrangian $\Lambda' \subset T^*\mathbb R^n$ described
in \eqref{eq:Ends}. Observe that the part of $\Lambda'$ lying over
$\{0\} \times\mathbb R ^p$ is the image of a Whitney type immersion
$S^p \to T^*\mathbb R^p \cong \{0\} \times T^*\mathbb R^p \subset
T^*\mathbb R^n$, obtained from the standard Whitney immersion
\begin{equation}
    \label{eq:Whitney-immersion}
    S^p \to T^*\mathbb R^p, \quad 
    (\x,y) \mapsto (\x, y\x) = (\x, \sqrt{1 - \vert \x \vert^2} \x)
\end{equation}
for $(\x,y) \in S^p \subset \mathbb R^p \times \mathbb R$, by
rescaling. This motivates the following definition:

\begin{defn}
    \label{defn:Whitney-deg}
    Let $L' \subset M$ be a Lagrangian submanifold containing the
    image of a Lagrangian immersion $\iota: D^{n-p} \times S^p \to M$
    which is an embedding away from $\{0\} \times S^p$ and such that
    $\check S = \iota(\{0\} \times S^p)$ has precisely one transverse
    double point. We call $\iota(D^{n-p} \times S^p) \subset L'$ a
    \emph{Whitney degeneration} if the following holds: There exists
    an embedded isotropic $p$-disc $\check D \subset M$ with boundary
    on $\check S$ and containing the double point of $\check S$ in its
    interior, together with a Weinstein neighbourhood
    $\mathcal N \cong (T\check D)^\omega/T \check D \oplus T^*\check
    D$ of $\check D$ such that upon a suitable symplectic
    identification of $\mathcal N$ with a subset of
    $T^*\mathbb R^{n-p} \times T^*\mathbb R^p$, $\check S$ is the
    image of a Whitney type immersion
    $S^p \to \{0\} \times T^*\mathbb R^p$ (see Figure
    \ref{fig:Whitney-sphere}).

\end{defn}

There a two basic symplectic invariants that one can associate to a
Whitney degeneration $\iota(D^{n-p} \times S^p) \subset L'$ with
$p > 0$. The first one is the symplectic area of the element of
$H_2(M,L')$ represented by the teardrop shown in Figure
\ref{fig:Whitney-sphere}. The second one is the pair of Maslov indices
of the discs created when resolving the double point by either of the
two topologically different ways of Lagrangian 0-surgery; we will
compute these Maslov indices for Whitney degenerations coming from
Lagrangian antisurgery in Section \ref{sec:monotonicity}. Note that,
in view of these computations, the second invariant yields an
obstruction to resolving a given Whitney degeneration
$\iota(D^{n-p} \times S^p) \subset L'$ by Lagrangian $p$-surgery.

As an example for the fact that containing a Whitney degeneration is
not sufficient for being able to perform Lagrangian $p$-surgery,
consider the Whitney sphere $S^n_{Wh} \subset T^*\mathbb R^n$
itself. It obviously contains Whitney degenerations
$\iota(D^{n-p} \times S^p)$ for every $0 \leq p \leq n-1$, but it is
not possible to perform Lagrangian $p$-surgery on $S^n_{Wh}$ for any
$p > 0$: This would lead to a compact embedded Lagrangian
$L \subset T^*\mathbb R^n$ with vanishing area class, which we know
not to exist. Indeed, $L$ would be diffeomorphic to
$S^{n-p-1} \times S^{p+1}$, which has $H_1(L) = 0$ for
$p \notin \{0,n-2\}$; in the case $p = n-2$, we would create an
isotropic 2-disc whose boundary generates $H_1(L) \cong \mathbb Z$,
and hence the area class would vanish as well.

\section{Desingularization}
\label{sec:singular-locus} 
Our aim in this section is to turn the immersed antisurgery cobordisms
$V: L' \leadsto L$ constructed in Section \ref{sec:defn-antisurgery} into
embedded cobordisms $V^\natural: L^\natural \leadsto L$ by
simultaneously resolving the singularities of $V$ and $L'$, and thus
prove Theorem \ref{main-thm:embedded-cob}.

\subsection{The singular loci of $\Gamma$ and $\Lambda'$}
\label{sec:singular-locus-Gamma-Lambda}
Recall from Section \ref{sec:immersed-cobordism} that the Lagrangian
$(k+1)$-handle $\Gamma: \Lambda' \leadsto \Lambda$, which serves as
the local model for the cobordisms corresponding to $k$-antisurgery,
is defined as the union of the graphs of $\pm dF$ for a function
$F: \mathbb R \times \mathbb R^n \to \mathbb R$. The locus $\Gamma^s$
where $\Gamma$ fails to be embedded is given by the points
$(x_0,\x) \in \mathrm{int}\,\mathfrak U$ where $dF(x_0,\x) = 0$, which
is where the graphs of $\pm dF$ intersect each other. Using this
description, we obtain the following proposition.

\begin{prop}
    \label{prop:sing-locus}
    The singular locus of the antisurgery handle $\Gamma: \Lambda'
    \leadsto \Lambda$ is given by
    \begin{equation}
        \label{eq:sing-locus}
        \Gamma^s = \{((x_0,0),(0,0)) \in T^*\mathbb R \times T^*\mathbb
        R^n~|~ x_0 \geq 1 - \delta\}.
    \end{equation}
\end{prop}
\begin{proof}
    Recall from Section \eqref{sec:immersed-cobordism} that
    \begin{equation*}
        dF = \frac{3}{2}f(x_0,\x)^{1/2}\big(\sigma'(x_0)\rho(r^2)dx_0 +
        (1 + \sigma(x_0)\rho'(r^2))dr^2 -ds^2\big).
    \end{equation*}
    Since $f(x_0, \x) > 0$ for
    $(x_0, \x) \in \mathrm{int} \,\mathfrak U$, the vanishing of
    $dF(x_0,\x)$ is equivalent to
    \begin{equation}
        \label{eq:dF=0}
        \begin{aligned}
            \sigma'(x_0)\rho(r^2)dx_0 &=0,\\
            (1 + \sigma(x_0)\rho'(r^2))dr^2 &=0,\\
            ds^2 & =0.
        \end{aligned}
    \end{equation}
    The conditions imposed on $\sigma$ and $\rho$ imply that these
    equations are simultaneously satisfied if and only if
    $x_0 \geq 1 - \delta$ and $\x = 0$, meaning that the singular
    locus of $\Gamma$ is as described in \eqref{eq:sing-locus}. 
    To see this, note that the conditions $\sigma(x_0) \geq 0$ and
    $-1/(1+\varepsilon) < \rho'(r^2) \leq 0$ imply
    $1 + \sigma(x_0) \rho'(r^2) > 0$, so the second equation in
    \eqref{eq:dF=0} can only hold when $dr^2 = 0$; together with the third
    equation $ds^2 = 0$, we conclude that $\x = 0$. Since $\rho(0) = 1$,
    the first equation simplifies to $\sigma'(x_0) dx_0 = 0$ and thus
    $x_0 \leq \delta$ or $x_0 \geq 1 - \delta$; since
    $(x_0,0) \notin \mathfrak U$ for $x_0 < 1 - \delta$, only the second
    of these possibilities leads to a solution of \eqref{eq:dF=0}.
\end{proof}

Recall that the part of $\Gamma$ that lies over $[1-\delta,\infty)
\times \mathbb R \subset T^*\mathbb R$ is cylindrical of the form
\begin{equation}
    \label{eq:15}
    \Gamma \vert_{[1-\delta,\infty) \times
        \mathbb R} = [1-\delta,\infty) \times \{0\} \times \Lambda' \subset T^*\mathbb R \times T^*\mathbb R^n.
\end{equation}
Proposition \ref{prop:sing-locus} therefore implies that the positive
end $\Lambda'$ of $\Gamma$ has a double point at $\x = 0$ and is
embedded away from that. The tangent spaces
\begin{equation}
    \label{eq:3}
    \lambda_\pm = T_{(0,0)}\Lambda_\pm'
\end{equation}
to the two sheets of $\Lambda'$ at the double point are spanned by
\begin{equation}
    \label{eq:tgt-space-0}
    \begin{aligned}
        \partial_{x_i} \pm 3 f_1(0)^{1/2}
        \partial_{y_i}, \quad & i = 1,\dots,k + 1,\\
        \partial_{x_i} \mp 3 f_1(0)^{1/2}
        \partial_{y_i}, \quad & i = k+2,\dots,n,\\
    \end{aligned}
\end{equation}
where $f_1(0) = f(1,0)$; since $f_1(0) \neq 0$, this shows that
$\lambda^+$ and $\lambda^-$ intersect transversely. 

The transverse double point of $\Lambda'$ can be removed by Lagrangian
0-surgery such as to produce an embedded Lagrangian submanifold
$\Lambda^\natural \subset T^*\mathbb R^n$. More interestingly, we will
see below that one can resolve the singular locus $\Gamma^s$ of $\Gamma$ and
turn the immersed cobordism $\Gamma: \Lambda' \leadsto \Lambda$ into
an \emph{embedded} Lagrangian cobordism $\Gamma^\natural:
\Lambda^\natural \leadsto \Lambda$.

\subsection{Implantation of a Lagrangian 1-handle}
\label{sec:model-case} To resolve the singular locus $\Gamma^s$ of
$\Gamma$, we replace a neighbourhood of $\Gamma^s$ with a Lagrangian
1-handle arising as a subset of the \emph{trace of surgery} cobordism
$\mathbb R^n \# i\mathbb R^n \leadsto (\mathbb R^n, i\mathbb R^n)$
constructed in \cite[Section 6.1]{Biran-Cornea--Lag-Cob-I}.\footnote{A
    Lagrangian 1-handle as needed here could also be obtained by
    modifying our construction of the 1-handle corresponding to
    $(n-1)$-antisurgery given in Section \ref{sec:immersed-cobordism}
    in such a way that the components of the positive end lie over
    disjoint curves in $T^*\mathbb R$, thus making the 1-handle
    embedded.
}
The explanations below and Proposition \ref{prop:trace-cob-revisited}
serve to justify that the cutting and isotoping in Biran--Cornea's
construction can be performed in such a way as to exactly match up the
resulting Lagrangian 1-handle with the given ``boundary condition''
that is created when one removes a neighbourhood of $\Gamma^s$ from
$\Gamma$.

\begin{figure}[t]
    \centering
    \includegraphics[scale=1]{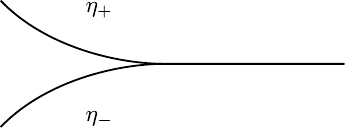}
    \caption{The singularity of $\Gamma$ is modelled by
        $\eta_+ \times \lambda_+ \cup \eta_- \times \lambda_-$, where
        $\eta_\pm \subset T^*\mathbb R$ are the curves depicted here,
        and $\lambda_\pm \subset T^*\mathbb R^n$ are transversely
        intersecting Lagrangian subspaces.}
    \label{fig:eta_pm}
\end{figure}

Let $\eta_\pm : \mathbb R \to T^*\mathbb R$ be curves given by
$\eta_\pm(x) = (x, \pm y(x))$, where
$y: \mathbb R \to \mathbb R_{\geq 0}$ is a smooth function such that
$y(x) > 0$ for $x < 0$ and $y(x) = 0$ for $x \geq 0$ (see Figure
\ref{fig:eta_pm}).
Let $\lambda_\pm \subset T^*\mathbb R^n$ be two transversely intersecting
Lagrangian subspaces. Then
\begin{equation*}
    W = \eta_+ \times \lambda_+ \cup  \eta_- \times \lambda_-
\end{equation*}
is an immersed Lagrangian submanifold of
$T^*\mathbb R \times T^*\mathbb R^n$ whose singular locus is
$\mathbb R_{\geq 0} \times \{0\}$. We view $W$ as an immersed
Lagrangian cobordism
$W: \lambda_- \cup \lambda_+ \leadsto (\lambda_-,\lambda_+)$ and will
use it as our local model for a neighbourhood of $\Gamma^s$.

\begin{prop}
    \label{prop:trace-cob-revisited}
    There exists an embedded Lagrangian cobordism $W^\natural:
    \lambda_- \# \lambda_+ \leadsto (\lambda_-, \lambda_+)$, which
    topologically is a 1-handle $\cong D^1 \times D^n$, such that
    $W^\natural$ and $W$ coincide outside of an arbitrarily small
    neighbourhood of the singular locus $\mathbb R_{\geq 0} \times
    \{0\}$ of $W$.
\end{prop}

\begin{proof} We will show how to perform the construction of
    Biran--Cornea's trace cobordism \cite{Biran-Cornea--Lag-Cob-I}
    corresponding to the Lagrangian surgery of $\lambda_\pm$ in such
    way that the result agrees with $W$ outside of an arbitrarily
    small neighbourhood of the singular locus
    $\mathbb R_{\geq 0} \times \{0\}$ of $W$. We assume that
    $\lambda_- = \mathbb R^n \times \{0\}$ and
    $\lambda_+ = \{0\} \times \mathbb R^n$ for notational simplicity.

    Choose a curve $\gamma : \mathbb R \to T^* \mathbb R$, $\gamma(t)
    = (a(t),b(t))$, as in Section \ref{sec:lagrangian-0-surgery} and
    let $L = \lambda_- \# \lambda_+ \subset T^*\mathbb R^n$ be the
    result of the corresponding 0-surgery of $\lambda_\pm$. Then define
    \begin{equation}
        \label{eq:model-cob-parametrization}
        \phi_\gamma: \mathbb R \times S^n \to T^*\mathbb R^{n+1}
    \end{equation}
    to be the composition of the map $\mathbb R \times S^n \to
    T^*\mathbb R^{n+1}$, $(t,\x) \mapsto (a(t) \x,b(t) \x)$, with a
    rotation of the first factor of $T^*\mathbb R^{n+1} = T^*\mathbb R
    \times T^* \mathbb R^n$ by $\frac{\pi}{4}$. Let
    \begin{equation*}
        W' = (\mathrm{im}\, \phi_\gamma) \vert_{\{(x, y) \in T^* \mathbb R
            | x \leq 0\}}
    \end{equation*}
    be the part of the image of $\phi_\gamma$ that lies over the
    half-plane $\{(x ,y) \in T^*\mathbb R ~|~ x \leq 0\}$ (using the
    notation introduced in Section
    \ref{sec:lagrangian-cobordisms}). Note that $W'$ is a manifold
    with boundary, and the boundary $\partial W'$ is given by
    \begin{equation*}
        L_0 = \{(0,0)\} \times L.
    \end{equation*}
    In the rest of the proof, we will describe how to adjust $W'$ such
    that a cylindrical end $\mathbb R_{\geq 0} \times L$ can be glued
    on, and such that the resulting Lagrangian looks like
    $W = \eta_+ \times \lambda_+ \cup \eta_- \times \lambda_-$ outside of
    a small neighbourhood of $\mathbb R_{\geq 0} \times \{0\} \subset
    T^* \mathbb R \times T^*\mathbb R^n$.

    To start, let $\mathcal N$ be a Weinstein neighbourhood of the
    Lagrangian $\mathbb R \times L \subset T^*\mathbb R^{n+1}$ which
    is of the form $\mathcal N = T^*\mathbb R \times \mathcal N_L$,
    where $\mathcal N_L \subset T^*\mathbb R^n$ is a Weinstein
    neighbourhood of $L$, and such that the map
    $\pi_{\mathcal N}: \mathcal N \to \mathbb R \times L$ induced by
    the canonical projection in the cotangent bundle is of the form
    \begin{equation}
        \label{eq:splitting}
        \pi_{\mathcal N} = \pi_{T^*\mathbb R} \times \pi_{\mathcal N_L},
    \end{equation}
    where $\pi_{T^*\mathbb R}$ and $\pi_{\mathcal N_L}$ are the
    corresponding maps for $\mathbb R \subset T^*\mathbb R$ and
    $L \subset \mathcal N$ (in particular,
    $\pi_{T^*\mathbb R}: T^*\mathbb R \to \mathbb R$ is simply the
    projection onto the first factor of $T^*\mathbb R \cong \mathbb R
    \times \mathbb R$).
    
    Let $U_0' \subset W'$ be a neighbourhood of $L_0 = \partial W'$ in
    $W'$. By shrinking it if necessary, we may assume that $U_0'$ lies
    entirely in the Weinstein neighbourhood $\mathcal N$ and that it
    is the graph of a closed 1-form $\alpha_0$ over the subset
    $U_0 = (-3 \varepsilon_0, 0] \times L \subset \mathbb R \times L$ for
    some small $\varepsilon_0 > 0$. Note that $\alpha_0$ is exact
    because its restriction to $L_0$ vanishes (as $L_0$ is contained
    in $\mathbb R \times L$) and because $U_0'$ retracts onto
    $L_0$. Moreover, $\alpha_0$ vanishes on $L_0$ and hence any
    primitive of $\alpha_0$ is constant on $L_0$. We denote by
    $g_0: U_0 \to \mathbb R$ the primitive of $\alpha_0$ which
    vanishes on $L_0$.

    Consider now the subset
    \begin{equation}
        \label{eq:components}
        U_0' \cap \big(T^*\mathbb R \times (B^n_{2 \kappa}(\lambda_+)^c \cup
        B^n_{2 \kappa}(\lambda_-)^c)\big)
    \end{equation}
    of $W'$, where $B^{n}_{2\kappa}(\lambda_\pm)^c$ denote the
    complements in $\lambda_\pm$ of the balls
    $B^{n}_{2 \kappa}(\lambda_\pm) \subset \lambda_\pm$ of radius
    $2 \kappa$ (it could also be written as
    $U_0' \cap \big(T^*\mathbb R^n \times B_{2 \kappa}(T^*\mathbb
    R^n)^c)$, where $B_{2 \kappa}(T^*\mathbb R^n)^c$ is the complement
    of a ball in $T^*\mathbb R^n$). This subset has two components
    which are contained in
    $\ell_\pm \times B^n_{2 \kappa}(\lambda_\pm)^c \subset T^*\mathbb
    R \times L$, where $\ell_\pm$ are the lines in $T^*\mathbb R$
    given by $y = \mp x$, i.e.\ the graphs of $d(\mp \frac{1}{2}x^2)$;
    in fact, these components are the graphs of $\alpha_0 = dg_0$ over
    $(-3 \varepsilon_0,0] \times B^n_{2 \kappa}(\lambda_\pm)^c \subset
    U_0$. Using this and the split nature $\eqref{eq:splitting}$ of
    $\pi_{\mathcal N}$, it follows that the restrictions of $g_0$ to
    $(-3 \varepsilon_0, 0] \times B^n_{2 \kappa}(\lambda_{\pm})^c$,
    which we denote by $g_0^\pm$, depend only on $x \in \mathbb R$ and
    are in fact given by $g_0^\pm(x) = \mp \frac{1}{2} x^2$, see
    Figure \ref{fig:cut-off}.

    Let now $\xi: U_0 \to \mathbb R$ be a cut-off function which
    depends only on $x \in \mathbb R$ and which satisfies
    $\xi(x) \equiv 0$ for
    $x \in (-3 \varepsilon_0, -2 \varepsilon_0]$, $\xi(x) \equiv 1$
    for $x \in [-\varepsilon_0, 0]$, and
    \begin{equation}
        \label{eq:condition-derivatives}
        \pm (\xi(x) g_0^\pm(x))' \leq \pm (g_0^\pm)'(x)
    \end{equation}
    for all $x$; it is not hard to verify that this last condition can
    be satisfied. Then denote by $\Phi: \mathcal N \to \mathcal N$ the
    time-one map of the Hamiltonian flow of $-\xi g_0$. We will use
    $\Phi$ to adjust $W'$ as required.
    \begin{figure}[t]
        \centering
        \includegraphics{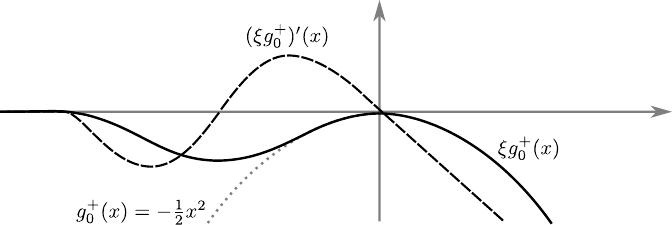}
        \caption{Construction of the function $\xi g_0$ whose
            Hamiltonian flow is used in the proof of Proposition
            \ref{prop:trace-cob-revisited} to adjust Biran--Cornea's
            Lagrangian 1-handle $W'$ to match the local model
            $W$ of the singular locus of the antisurgery cobordism
            away from a small neighborhood of its singular locus.
        }
        \label{fig:cut-off}
    \end{figure}
    
    By construction, $\Phi$ takes
    $U_0' \cap ((-\varepsilon_0, 0] \times T^*\mathbb R^n)$ to
    $(-\varepsilon_0, 0] \times L$ and leaves
    $U_0' \cap ((-3 \varepsilon_0, -2 \varepsilon_0] \times T^*\mathbb
    R^n)$ fixed. We can therefore extend $\Phi$ to a map
    $\hat \Phi: W' \to T^*\mathbb R \times T^*\mathbb R^n$ defined on all
    of $W'$ that leaves
    $W' \cap ((- \infty, -2 \varepsilon_0] \times T^*\mathbb R^n)$
    fixed and that takes
    $W' \cap ((-\varepsilon_0, 0] \times T^*\mathbb R^n)$ to
    $\mathbb R \times L$. The last statement implies that we can extend
    $\hat \Phi(W')$ by a cylindrical end such that it becomes a valid
    Lagrangian cobordism.
    
    Note that the restrictions of $\Phi$ to the components of
    \eqref{eq:components} are of the form
    $\phi_\pm \times \mathrm{id}$, where $\phi_\pm$ are the time-one
    maps of the Hamiltonian flows of the restrictions of $\xi g_0$ to
    the corresponding components of
    $U_0 \cap (T^*\mathbb R \times (B^{n}_{2 \kappa}(\lambda_+)^c \cup
    B^{n}_{2 \kappa}(\lambda_-)^c))$. In other words, $\Phi$ moves
    these sets in the direction of the fibres of $T^*\mathbb R$. The
    extended map $\hat \Phi$ therefore takes
    $W' \cap (T^*\mathbb R \times (B^{n}_{2 \kappa}(\lambda_+)^c \cup
    B^{n}_{2 \kappa}(\lambda_-)^c))$ to
    $\tilde \eta_+ \times B^{n}_{2 \kappa}(\lambda_+)^c \cup
    \tilde \eta_- \times  B^{n}_{2 \kappa}(\lambda_-)^c$, where
    $\tilde \eta_\pm = \phi_\pm(\eta_\pm)$. It follows from the
    inequalities \eqref{eq:condition-derivatives} that the curves
    $\tilde \eta_\pm$ are entirely contained in the upper resp.\ lower
    half-planes and only intersect along the $x$-axis; that is, up to
    a horizontal shift they are of the type describe right before this
    proposition. Using an appropriate symplectomorphism
    $\psi: T^*\mathbb R \to T^*\mathbb R$ of the form
    $(x, y) \mapsto (f(x), g(x,y))$, we can match up the curves
    $\tilde \eta_\pm$ with $\eta_\pm$. (Such a symplectomorphism
    exists provided that the functions $y, \tilde y$ defining
    $\eta_\pm, \tilde \eta_\pm$ via $\eta_\pm(x) = (x, \pm y(x))$
    etc.\ satisfy $\lim_{x \to 0} y(x)/\tilde y(x) \to 1$, see Lemma
    \ref{lem:matching-curves}; this condition can be guaranteed by
    choosing the cut-off function $\xi$ appropriately.)

    \begin{figure}[t]
        \centering
        \includegraphics[scale=1]{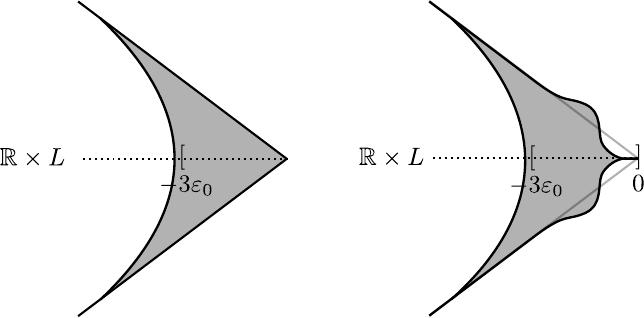}
        \caption{Projection of $W'$ resp.\ $\Phi(W')$ to
            $T^* \mathbb R$, where $W'$ is Biran--Cornea's Lagrangian
            1-handle and $\Phi$ is the time-one map of the Hamiltonian
            flow of the function $\xi g_0$ constructed in the proof of
            Proposition \ref{prop:trace-cob-revisited} in order to
            adjust $W'$ as required.}
        \label{fig:projection-adjusted-cob}
    \end{figure}

    Applying $\psi \times \textrm{id}$ to $\hat \Phi(W')$ and gluing on a
    cylindrical end hence yields a Lagrangian cobordism $W^\natural$ that
    agrees with $W$ outside of an arbitrarily small neighbourhood of the
    singular locus of $W$.
\end{proof}

\begin{lem}
    \label{lem:matching-curves}
    Let $\eta^i: \mathbb R \to T^*\mathbb R$, $x \mapsto (x, y^i(x))$,
    $i = 0,1$, be curves as described above for which $\lim_{x \to 0-}
    (y^0(x)/y^1(x)) = 1$. Then there exists a symplectomorphism
    $\Psi: T^*\mathbb R \to T^*\mathbb R$  taking $\eta^0_\pm$ to $\eta^1_\pm$.
\end{lem}
\begin{proof}
    We will construct a diffeomorphism of the form
    $\Psi(x,y) = (f(x), g(x,y))$. For such $\Psi$ to be a
    symplectomorphism, we need
    \begin{equation*}
        \left \vert
            \begin{pmatrix}
                \partial_x f & 0\\
                \partial_x g & \partial_y g
            \end{pmatrix}
        \right \vert
        = \partial_x f \cdot \partial_yg = 1
    \end{equation*}
    everywhere, which implies
    \begin{equation*}
        g(x,y) = \frac{y}{\partial_x f(x)}
    \end{equation*}
    up to an additive constant. In order for $\Psi$ to take
    $\pm \eta^0$ to $\pm \eta^1$, we need
    \begin{equation*}
        g(x, \pm y^0(x)) = \pm y^1(f(x)),
    \end{equation*}
    for $x < 0$ and $f(x) = x$ for $x \geq 0$. In view of the previous
    equation, this is equivalent to the ordinary differential equation
    \begin{equation*}
        \partial_x f(x) = \frac{y^0(x)}{y^1(f(x))}
    \end{equation*}
    for $x < 0$. Provided that $\lim_{x \to 0-} (y^0(x)/y^1(x)) = 1$,
    this ODE has a solution which extends to a function
    $f: \mathbb R \to \mathbb R$ with $f(x) = x$ for $x \geq 0$.
\end{proof}

\subsection{Surgery of the singular locus of $\Gamma$}
\label{sec:conclusion}
In suitable Darboux coordinates, a neighbourhood
$\mathcal U(\Gamma^s)$ of the singular locus $\Gamma^s$ of $\Gamma$
looks like
$(\eta_- \times B^{n}_{r}(\lambda_-)) \cup (\eta_+ \times
B^n_{r}(\lambda_+))$, where $B^n_{r}(\lambda_\pm)$ denotes a ball of
radius $r$ in $\lambda_\pm = T_{(0,0)}\Lambda'_\pm$ and with curves
$\eta_\pm$ as described in the previous subsection, up to a shift to
the right by $1 - \delta$ and up to restricting them from $\mathbb R$
to $[-\nu, \infty)$ for some small $\nu > 0$. Explicitly, we can take
$\eta_\pm(x) = (x, \pm \frac{3}{2} (\sigma(x) - 1)^{1/2} \sigma'(x))$.

Having set up such an identification, we can use Proposition
\ref{prop:trace-cob-revisited} to replace this neighbourhood with a
corresponding piece of the 1-handle
$W^\natural: \lambda_- \# \lambda_+ \leadsto (\lambda_-,
\lambda_+)$. More precisely, we can take a 1-handle $W^\natural$ as
described there with $2 \kappa < r$ and replace
\begin{equation}
    \label{eq:cut-paste-1}
    \mathcal U(\Gamma^s) \cap \Big(T^*\mathbb R \times \big(B^{n}_{2
        \kappa}(\lambda_+)^c \cup B^n_{2 \kappa}(\lambda_-)^c\big)\Big)^c
\end{equation}
with
\begin{equation}
    \label{eq:cut-paste-2}
    W^\natural \cap \Big(T^* \mathbb R \times \big(B^{n}_{2
        \kappa}(\lambda_+)^c \cup B^n_{2
        \kappa}(\lambda_-)^c\big)\Big)^c.
\end{equation}
We assume here tacitly that the 1-handle $W^\natural$ has been pruned
analogously to restricting the curves $\eta_\pm$ from $\mathbb R$ to
$[-\nu, \infty)$, i.e., by removing the parts lying over
$\eta_\pm((-\infty, \nu))$. With this in mind, the boundaries of
\eqref{eq:cut-paste-1} and \eqref{eq:cut-paste-2} are both given by
$\partial^v \cup \partial^h$ with
\begin{align}
  \label{eq:attaching-region}
  \begin{split}
      \partial^v &= \{\eta_+(-\nu)\} \times B^n_{2 \kappa}(\lambda_+) \cup
      \{\eta_-(-\nu)\}
      \times B^n_{2 \kappa}(\lambda_-)\\
      \partial^h &= \eta_+([-\nu, \infty)) \times \partial B^n_{2
          \kappa}(\lambda_+) \cup \eta_-([-\nu, \infty))
      \times \partial B^n_{2 \kappa}(\lambda_-),
  \end{split}
\end{align}
which makes the cutting and pasting along this ``attaching region'' possible.

The outcome of this operation is an embedded Lagrangian cobordism
\begin{equation}
    \label{eq:model-cob}
    \Gamma^\natural: \Lambda^\natural \leadsto \Lambda,
\end{equation}
where $\Lambda^\natural$ is the result of resolving the double point
of $\Lambda'$ by an application of Lagrangian 0-surgery which is
locally modelled on $\lambda_- \# \lambda_+$. By choosing the curve
$\gamma$ used in the proof of Proposition
\ref{prop:trace-cob-revisited} sufficiently small, we can guarantee
that $\Gamma^\natural$ and $\Gamma$ coincide outside of an arbitrarily
small neighbourhood of $\Gamma^s$ in $T^*\mathbb R^{n+1}$.

\subsection{Proof of Theorem \ref{main-thm:embedded-cob}}
\label{sec:proof-theorem}
Suppose that we have a Lagrangian $L \subset M$ with an embedded
sphere $S \subset L$ and an isotropic surgery disc $D$ for
$S$. Repeating the construction in Section \ref{sec:defn-antisurgery},
but replacing $\Gamma$ by $\Gamma^\natural$, the desingularized
Lagrangian cobordism constructed in Section \ref{sec:conclusion}, and
consequently $\Lambda$ by $\Lambda^\natural$, we produce a Lagrangian
cobordism
\begin{equation*}
    V^\natural: L^\natural \leadsto L,
\end{equation*}
where $L^\natural$ is a Lagrangian obtained by resolving the double
point created when performing antisurgery on $L$ along the isotropic
disc $D$ by an application of Lagrangian 0-surgery.
The cobordism $V^\natural$ is embedded if $L$ is embedded, and it can
be arranged to agree with with the corresponding $V: L' \leadsto L$
outside of an arbitrarily small neighbourhood of its singular locus by
choosing the parameter $\kappa$ sufficiently small.
Topologically, this cobordism is the concatenation of the traces
corresponding to first performing $k$-surgery on $L$ and then
$0$-surgery on the result $L'$ of the first step; in other words,
$V^\natural$ is obtained from $[0,1] \times L$ by first attaching a
$(k+1)$-handle and then a $1$-handle.

This completes the proof of Theorem
\ref{main-thm:embedded-cob}. \qed

\subsection{The Lagrangian isotopy class of the resolution}
\label{sec:Lagr-isotopy-class}
Recall from Section \ref{sec:lagrangian-0-surgery} that there are
\emph{two} families of resolutions of the double point of $\Lambda'$
by Lagrangian $0$-surgery which depend on a choice of order of sheets
near the double point (and correspond to Lagrangian isotopy classes of
0-surgery models). It is important to note that in the construction
leading to the desingularized antisurgery cobordism
$V^\natural: L^\natural \leadsto L$ whose existence is asserted by
Theorem \ref{main-thm:embedded-cob} we do \emph{not} have a choice
regarding which of these families the desingularized end
$\Lambda^\natural$ belongs to, as Proposition
\ref{prop:trace-cob-revisited} only gives the existence of a
Lagrangian cobordism
$(\lambda_-,\lambda_+) \leadsto \lambda_- \# \lambda_+$ (as opposed to
$(\lambda_-,\lambda_+) \leadsto \lambda_+ \# \lambda_-$).

Recall from Section \ref{sec:lagrangian-0-surgery-1} that in the case
$k = n-1$, both $L$ and $L^\natural$ are resolutions of $p \in L'$,
the double point created when applying $(n-1)$-antisurgery to $L$, by
Lagrangian 0-surgery. In this situation, we have:

\begin{prop}
    \label{prop:distinct_0_surgeries}
    The ends of the desingularized cobordism
    $V^\natural: L^\natural \leadsto L$ resulting from
    $(n-1)$-antisurgery on $L$ belong to \emph{distinct} families of
    resolutions of $p \in L'$ by Lagrangian 0-surgery.
\end{prop}
\begin{proof}
    As observed in Section \ref{sec:lagrangian-0-surgery-1}, $\Lambda$
    belongs to the family of resolutions of $\Lambda'$ locally given
    by $\lambda_+ \# \lambda_-$, while $L^\natural$ belongs to the
    family of resolutions locally given by $\lambda_- \#
    \lambda_+$.
\end{proof}

\subsection{The size of the resolution.}
\label{sec:upper-bound-resolution}
Let $L'$ be the singular end of a Lagrangian cobordism
$V: L' \leadsto L$ arising from antisurgery with parameter
$\varepsilon$.  In the following we present evidence that the size of
a Lagrangian 0-surgery (in the sense of Definition
\ref{defn:size-of-0-surgery}) which can be applied to $L$ in such a
way that the cobordism can be desingularized simultaneously, is upper
bounded by $2 \varepsilon^{3/2}$.

\begin{figure}[t]
    \centering
    \includegraphics{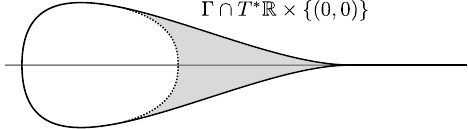}
    \caption{The intersection of the antisurgery cobordism $\Gamma$
        with the plane
        $T^* \mathbb R \times \{(0,0)\} \subset T^* \mathbb R \times
        T^* \mathbb R^n$. The dotted line shows the intersection of
        the Lagrangian 1-handle used for desingularization with the plane
        $T^*\mathbb R \times \{(0, 0)\}$.
        }
    \label{fig:teardrop}
\end{figure}

Consider the intersection of the model cobordism $\Gamma$ with the
plane
$T^*\mathbb R \times \{(0, 0)\} \subset T^*\mathbb R \times T^*
\mathbb R^n$.  By setting $\x = (0, 0) \in T^*\mathbb R^n$ in
\eqref{eq:dF}, one can see that this intersection consists of a curve
that bounds a teardrop-shaped region as shown in Figure
\ref{fig:teardrop}. The Lagrangian 1-handle that we use in order to
resolve the singularity of $\Gamma$ also intersects
$T^*\mathbb R^n \times \{(0, 0)\}$ in a curve, which is obtained by
restricting the map
$\phi_\gamma: \mathbb R \times S^n \to T^*\mathbb R^{n+1}$
\eqref{eq:model-cob-parametrization} modelling this handle to
$I \times \{(1, 0, \dots, 0)\} \subset \mathbb R \times S^n$ for a
suitable interval $I \supset [-\kappa, \kappa]$ (with $\kappa$ as in
the description of the local model for Lagrangian 0-surgery
\eqref{eq:surgery-curve-spec}); after the implanting the 1-handle,
this curve lies inside the teardrop (the part of the 1-handle that
projects to the coordinate axes in the local model projects to the
boundary curve of the teardrop after implantation).

This picture suggest that the area of the teardrop is an upper bound
for the size of the Lagrangian 0-surgery (the area of the shaded
region in Figure \ref{fig:teardrop}) that we can perform on the
positive end of the antisurgery cobordism in such a way that we can
simultaneously desingularize the latter. A simple computation using
again \eqref{eq:dF} shows that the area of the teardrop is
$2 \varepsilon^{3/2}$.

\subsection{Orientability}
\label{sec:orientability}
The result of abstract $k$-surgery on an orientable manifold $L$ is
always orientable if $k \geq 1$, since the $D^{k+1} \times S^{n-k-1}$
we glue in has a connected boundary in that case (or rather, every
component has a connected boundary---this includes the case $k =
n-1$). In the case $k = 0$, orientability depends on whether $D^1
\times S^{n-1}$ is glued consistently along its two boundary
components. Let
\begin{equation*}
        P^n = S^1 \times S^{n-1} \quad \text{and} \quad
        Q^n = D^1 \times S^{n-1}/\sim\,,
\end{equation*}
where $\sim$ identifies $\{1\} \times S^{n-1}$ with
$\{-1\} \times S^{n-1}$ using an orientation-reversing involution of
$S^{n-1}$ such as
$(x_1, x_2, \dots, x_n) \mapsto (-x_1, x_2, \dots, x_n)$ (so $Q^n$ is
the mapping torus of such an involution). The result of $0$-surgery on
$L$ is diffeomorphic to $L \# P^n$ in the orientable case and to
$L \# Q^n$ in the non-orientable case.

Returning to the Lagrangian setting, assume that the $L$ we start with
is orientable. When passing from $L$ to $L'$ by $k$-antisurgery, we
replace an embedded copy of $S^k \times D^{n-k}$ (a subset of
$\Lambda$) by an immersed copy of $D^{k+1} \times S^{n-k-1}$ with a
transverse double point (a subset of $\Lambda'$), and we resolve this
double point by 0-surgery when passing from $L'$ to $L^\natural$. Thus
$L^\natural$ is obtained from $L$ by replacing an embedded copy of
$S^k \times D^{n-k-1}$ by an embedded copy (a subset of
$\Lambda^\natural$) of either
\begin{equation*}
    (D^{k+1} \times S^{n-k-1}) \# P^n \quad \text{or} \quad (D^{k+1} \times S^{n-k-1}) \# Q^n
\end{equation*}
in the case $k < n-1$ (for $k=n-1$,
$D^{k+1} \times S^{n-k-1} = D^n \times S^0$ is not connected). The
first possibility leads to $L^\natural$ being orientable, the second
to $L^\natural$ being non-orientable; the next proposition states when
which of these alternatives holds.

\begin{prop}
    \label{prop:orientability}
    Let $L$ be an orientable Lagrangian and let $L^\natural$ be the
    result of $k$-antisurgery on $L$ and subsequent desingularization,
    for some $k$ with $0 < k < n-1$. Then $L^\natural$ is orientable
    if $k$ is odd and non-orientable if $k$ is even.
\end{prop}

\begin{proof}
    Note first that the result of abstract $k$-surgery on $L$ is
    always orientable when $0 < k < n-1$. 
    Whether $L^\natural$ is orientable or not is a local question that
    depends on whether the copy of $\mathbb R \times S^{n-1}$ we glue
    in when performing Lagrangian 0-surgery on $\Lambda'$ matches up
    orientations of the sheets $\Lambda_\pm'$ which induce the same
    orientation of $\Lambda'$. Assuming that such orientations for
    $\Lambda_\pm'$ have been chosen, it follows from \cite[Theorem
    4]{Polterovich--Surgery-of-Lagrangians}\footnote{When applying
        Polterovich's theorem, one needs to take into account that our
        convention for the direction of surgery handles is opposite to
        the one in \cite{Polterovich--Surgery-of-Lagrangians};
        according to which our handle is directed ``from $\lambda_+$
        to $\lambda_-$''; this becomes clear by inspection of the
        description of the model handle in \cite[Section
        1.8]{Polterovich--Surgery-of-Lagrangians}.}  that the result
    is orientable iff
    $(-1)^{n(n-1)/2 + 1}\Lambda'_+\cdot \Lambda'_- = 1$, where
    $\Lambda'_+ \cdot \Lambda'_-$ denotes the intersection index with
    respect to the symplectic orientation $\mathfrak o_{\omega}$ of
    the ambient manifold.
    
    To find such orientations for $\Lambda'_\pm$, consider first the
    orientations $\mathfrak o_\pm$ which the projections
    $\Lambda_\pm' \to \mathbb R^n \times \{0\}$ to the zero-section
    match up with the standard orientation of
    $\mathbb R^n \times \{0\}$. These are represented by the ordered
    bases of $\lambda_\pm = T_{(0,0)}\Lambda_\pm'$ listed in
    \eqref{eq:tgt-space-0}. By continuously deforming these to bases
    of $T_{(\x, \y)} \Lambda'$ with $(\x, \y)$ in the boundary of
    $\Lambda_\pm'$ (such as
    $(\x, \y) = ((0,\dots,0, \sqrt{\varepsilon}),(0,\dots,0))$), one
    can see that $\mathfrak o_\pm$ induce \emph{different}
    orientations of $\Lambda'$; a choice of orientations inducing the
    \emph{same} orientation of $\Lambda'$ is hence given by
    $\mathfrak o_+$ and $- \mathfrak o_-$.

    To compute $\Lambda_+'\cdot \Lambda_-'$ when $\Lambda_\pm'$ carry
    the orientations $\pm \mathfrak o_\pm$, note first that
    $\mathfrak o_+ \oplus \mathfrak o_-$ is represented by the ordered
    basis 
    \begin{align*}
      \big(&\partial_{x_1} +
             \partial_{y_1},\dots, \partial_{x_{k+1}} +
             \partial_{y_{k+1}},
             \partial_{x_{k+2}} -
             \partial_{y_{k+2}},\dots,\partial_{x_n}
             - \partial_{y_n},\\
      &\partial_{x_1} -
             \partial_{y_1},\dots, \partial_{x_{k+1}} -
             \partial_{y_{k+1}},
             \partial_{x_{k+2}} +
             \partial_{y_{k+2}},\dots,\partial_{x_n}
             + \partial_{y_n} \big)
    \end{align*}
    of $\mathbb R^{2n} \cong T_{(0,0)}T^*\mathbb R{^n}$. The
    linear map taking this basis to the symplectic basis
    \begin{equation*}
        \big(\partial_{x_1}, \partial_{y_1},\dots, \partial_{x_n},
        \partial_{y_n}\big)
    \end{equation*}
    has determinant
    $(-1)^{n(n-1)/2 + k + 1} 2^n$, and thus
    $\mathfrak o_+ \oplus \mathfrak o_- = (-1)^{n(n-1)/2 + k + 1}
    \mathfrak o_\omega$. It follows that
    \begin{equation*}
        \mathfrak o_+ \oplus (-\mathfrak o_-) = (-1)^{n(n-1)/2 + k}
        \mathfrak o_\omega,
    \end{equation*}
    and therefore
    \begin{equation*}
        (-1)^{n(n-1)/2 + 1} \Lambda_+' \cdot \Lambda_-' = (-1)^{k + 1},
    \end{equation*}
    from which the claimed statement follows.
\end{proof}

\begin{rmk}
    For $k = 0$, it is possible that the manifold $L'$ resulting from
    abstract $0$-surgery on an orientable $L$ is already
    non-orientable, in which case $L^\natural$ is also non-orientable;
    if, however, $L'$ is orientable, the statement of Proposition
    \ref{prop:orientability} applies to the resulting
    $L^\natural$. For $k = n-1$, it is possible that
    $L \smallsetminus (S^{n-1} \times D^1)$ is disconnected, in which
    case $L^\natural$ is orientable.
\end{rmk}

\subsection{Computation of Maslov indices}
\label{sec:monotonicity}
Consider again $\Lambda^\natural$, i.e., the resolution of the
singular end $\Lambda'$ of the $(k+1)$-handle
$\Gamma: \Gamma' \leadsto \Gamma$ corresponding to $k$-antisurgery for
which a desingularized cobordism
$\Gamma^\natural: \Lambda^\natural \leadsto \Lambda$
\eqref{eq:model-cob} exists. Note that $\Lambda^\natural$ is
diffeomorphic to either $(D^{k+1} \times S^{n-k-1}) \# P^n$ or
$(D^{k+1} \times S^{n-k-1}) \# Q^n$, with $P^n$ and $Q^n$ as defined
in Section \ref{sec:orientability}. For $k$ satisfying
$0 \leq k < n-2$, we have
\begin{equation*}
    H_2(T^*\mathbb R^n, \Lambda^\natural) \cong \mathbb Z
\end{equation*}
and there exists a preferred generator
$\sigma \in H_2(T^*\mathbb R^n, \Lambda^\natural)$ characterized by
the positivity of its symplectic area. In the proof Proposition
\ref{prop:Maslov-indices} below we will describe $\sigma$ by
describing a loop $\ell$ on $\Lambda^\natural$ representing
$\partial \sigma \in H_1(\Lambda^\natural) \cong H_2(T^* \mathbb R^n,
\Lambda^\natural)$ (this description also extends to the case
$k=n-2$).

Proposition \ref{prop:Maslov-indices} computes the Maslov index
$\mu(\sigma)$ of this generator, and in particular shows that
$\mu(\sigma)$ is non-positive for every $k$ satisfying
$1 \leq k \leq n-2$. As a consequence, Lagrangians resulting from
$k$-antisurgery and subsequent desingularization for $k$ in that range
are never monotone.

\begin{prop}
    \label{prop:Maslov-indices}
    Assume that $0 \leq k < n-1$. The Maslov index of
    $\sigma \in H_2(T^*\mathbb R^n, \Lambda^\natural)$ is given by
    $\mu(\sigma)= 1 -k$.
\end{prop}

    \begin{figure}[t]
        \centering
        \includegraphics[scale=1]{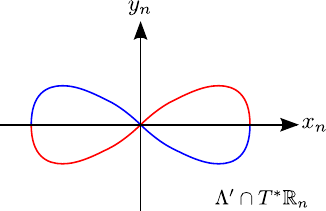}\qquad
        \includegraphics[scale=1]{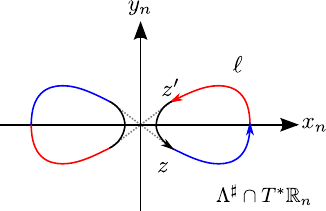}
        \caption{The intersections of the singular end $\Lambda'$ of
            the antisurgery cobordism and its resolution
            $\Lambda^\natural$ with the plane
            $T^* \mathbb R_n \subset T^* \mathbb R^n$. The loop $\ell$
            represents the boundary of the generator
            $\sigma \in H_2(T^*\mathbb R^n, \Lambda^\natural)$ whose
            Maslov index is computed in Proposition
            \ref{prop:Maslov-indices}.}
        \label{fig:intersection-surgered}
    \end{figure}

\begin{proof}
    Denote by $\mathbb R_i$ the $x_i$-coordinate subspace of
    $\mathbb R^n$ and by $T^*\mathbb R_i$ the $(x_i,y_i)$-coordinate
    subspace of $T^*\mathbb R^n$, for $i = 1,\dots,n$. To compute
    $\mu(\sigma)$, we will represent $\partial \sigma$ by a loop
    $\ell$ in $\Lambda^\natural \cap T^*\mathbb R_n$ and compute how
    the tangent spaces to $\Lambda^\natural$ twist as we traverse
    $\ell$.

    Recall from Section \ref{sec:immersed-cobordism} that we described
    the handle $\Gamma: \Lambda' \leadsto \Lambda$ as the union of the
    graphs of $\pm dF$ for a certain function $F$. Specializing the
    formula \eqref{eq:dF} for $dF$ to the case $x_0 = 1$, we see that
    the differential of $F_1 = F(1, \cdot)$ is
    \begin{equation}
        \label{eq:dF_1}
        dF_1 = \frac{3}{2}f_1(\x)^{1/2}\big((1 + (1 + \varepsilon)\rho'(r^2))dr^2 -ds^2\big)
    \end{equation}
    with
    $f_1(\x) = f(1, \x) = r^2 + (1 + \varepsilon)\rho(r^2) - s^2 - 1$
    (see Section \ref{sec:constr-cobord}). In particular, at points of
    the form $\x = (0,\dots,0,x_n) \in \mathbb R_n$ is given by
    $dF_1(\x) = -3 (\varepsilon - x_n^2)^{1/2} x_n dx_n$, and hence
    \begin{equation*}
        \Lambda' \cap T^* \mathbb R_n = \{(x_n,\mp 3(\varepsilon -
        x_n^2)^{1/2}x_n) \in T^*\mathbb R_n~|~ x_n^2 \leq \varepsilon\},
    \end{equation*}
    as depicted in the left part of Figure
    \ref{fig:intersection-surgered}, where the blue segment
    corresponds to $+dF_1$ and the red segment to
    $-dF_1$. Differentiation of \eqref{eq:dF_1} shows that the tangent
    space to $\Lambda'_\pm = \mathrm{graph}(\pm dF_1)$ over
    $\x = (0,\dots,0,x_n) \in \mathbb R_n$ is spanned by
    \begin{equation}
        \label{eq:10}
        \begin{aligned}
            \partial_{x_i} \pm 3 f_1(\x)^{1/2}
            \partial_{y_i}, \quad & i = 1,\dots,k+1\\
            \partial_{x_i} \mp 3 f_1(\x)^{1/2}
            \partial _{y_i}, \quad &i = k+2,\dots,n-1\\
            \partial _{x_n} \mp 3 \left( f_1(\x)^{1/2} -
                f_1(\x)^{-1/2}x_n^2 \right) \partial_{y_n}, \quad & i
            = n.
        \end{aligned}
    \end{equation}
    Note that the last vector is proportional to
    $f_1(\x)^{1/2}\partial_{x_n} \mp 3 \left( f_1(\x) - x_n^2
    \right)\partial_{y_n}$, so it approaches a multiple of
    $\partial_{y_n}$ as $x_n^2 \to \varepsilon $ (which is the
    conormal direction to $\{f_1(\x) = 0\} \subset \mathbb R^n$ at
    $\x = (0,\dots,0,\pm \varepsilon^{1/2})$). Moreover, \eqref{eq:10}
    shows that all these tangent spaces split as direct sums of
    1-dimensional subspaces of the 2-dimensional planes
    $T^*\mathbb R_i \subset T^*\mathbb R^n$; in particular, the
    tangent spaces to $\Lambda'_\pm$ at the origin are of the form
    $\lambda_\pm = \lambda_\pm^1 \times \dots \times \lambda_\pm^n$,
    where $\lambda_\pm^i$ is a 1-dimensional subspaces of
    $T^*\mathbb R_i$ for $i = 1, \dots, n$.
    
    To see how the 0-surgery by which we pass from $\Lambda'$ to
    $\Lambda^\natural$ affects the picture, recall from Section
    \ref{sec:lagrangian-0-surgery} that in order to obtain
    $\Lambda^\natural$, we use a symplectomorphism
    $\Phi: T^*\mathbb R^n \to T^*\mathbb R^n$ which identifies
    neighbourhoods of the origin in $\mathbb R^n \times \{0\}$ resp.\
    $\{0 \} \times \mathbb R^n$ with neighbourhoods of the origin in
    $\Lambda'_- = \{(\x, -dF_1(\x)) \in T^*\mathbb R^n~|~ \x \in
    \mathfrak U_1\}$ resp.\
    $\Lambda'_+ = \{(\x,+ dF_1(\x)) \in T^*\mathbb R^n~|~ \x \in
    \mathfrak U_1\}$ (cf.\ the explanation in Section
    \ref{sec:Lagr-isotopy-class}). After perturbing $\Lambda'$ a bit
    such that it agrees with $\lambda_- \cup \lambda_+$ near the
    origin, we may assume that $\Phi$ is linear and of the form
    $\phi^{1} \times \dots \times \phi^n$ with linear maps
    $\phi^i: T^*\mathbb R \to T^*\mathbb R_i$ that take
    $\mathbb R \times \{0\}$ to $\lambda_-^i$ and
    $\{0\} \times \mathbb R$ to $\lambda_+^i$.  Using such an
    identification, we glue in a Lagrangian copy of
    $D^1 \times S^{n-1}$ which is given by a part of the image of the
    map $h_\gamma: \mathbb R \times S^{n-1} \to T^*\mathbb R^n$
    defined in \eqref{eq:surgery-profile}. Figure
    \ref{fig:intersection-surgered} shows what
    $\Lambda^\natural \cap T^*\mathbb R_n$ looks like.

    We now compute $\mu(\sigma)$. Let $\ell$ be the loop on
    $\Lambda^\natural \cap T^*\mathbb R_n$ shown in Figure
    \ref{fig:intersection-surgered} traversed in counterclockwise
    direction, which is a representative of
    $\partial \sigma \in H_1(\Lambda^\natural)$. The black segment of
    $\ell$ is the image of the restriction of
    $\Phi \circ h_\gamma: I \times S^{n-1} \to T^*\mathbb R^n$ to
    $I \times \{(0, \dots, 0, -1)\}$, where $I \subset \mathbb R$ is a
    small interval containing 0. 
    Differentiation of \eqref{eq:surgery-profile} shows that the
    tangent space of $\Lambda^\natural$ at
    $\Phi \circ h_\gamma(t,(0,\dots,0,-1))=
    \Phi(-a(t)e_n, -b(t)e_n)$ is spanned by
    \begin{equation}
        \label{eq:13}
        \begin{aligned}
            \Phi \circ Dh_\gamma(\partial_{x_i}) &= \phi^i(a(t)e_i,
            b(t)e_i), \quad i = 1,\dots,n-1, \\
            \Phi \circ Dh_\gamma(\partial_t) &= \phi^n(-\dot a(t) e_n,
            - \dot b(t) e_n);
        \end{aligned}
    \end{equation}
    here $e_1, \dots, e_n$ are the standard basis vectors of
    $\mathbb R^n$.

    The formulas \eqref{eq:10} and \eqref{eq:13} show that the loop in
    the Lagrangian Grassmannian $Gr_L(T^*\mathbb R^n)$ induced by
    $\ell$ is contained in
    $Gr_L(T^*\mathbb R_1) \times \cdots \times Gr_L(T^*\mathbb R_n)
    \subset Gr_L(T^*\mathbb R^n)$, i.e., it is a direct sum of loops
    in $Gr_L(T^*\mathbb R_j)$.  Figure \ref{fig:Maslov-} indicates
    what the pieces of these loops corresponding to the pieces of
    $\ell$ look like; the pictures can be deduced from formulas
    \eqref{eq:10}, \eqref{eq:13}. 
    One can read off from these pictures the Maslov indices of the
    loops in $Gr_L(T^*\mathbb R_j)$, which for $j = 1, \dots, k+1$ are
    $-1$, for $j = k+2, \dots, n-1$ are 0, and for $j = n$ is 2; the
    Maslov index of $\sigma$ is their sum, $\mu(\sigma) = 1 - k$.
    \begin{figure}[t]
        \centering    
        \includegraphics[scale=.8]{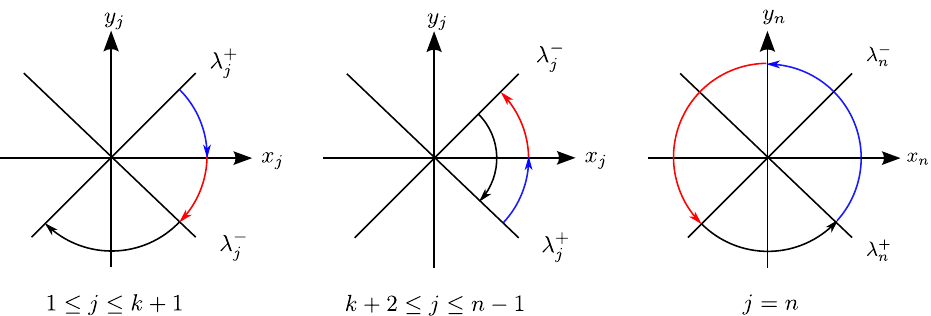}
        \caption{Computation of the Maslov index $\mu(\sigma)$ of the
            generator
            $\sigma \in H_2(T^* \mathbb R^n, \Lambda^\natural)$ whose
            boundary is represented by the loop $\ell$ in Figure
            \ref{fig:intersection-surgered}. The corresponding loop in
            $Gr_L(T^* \mathbb R^n)$ is a direct sum of loops in
            $Gr_L(T^*\mathbb R_j)$, $j = 1, \dots, n$, which are shown
            here.}
        \label{fig:Maslov-}
    \end{figure}
\end{proof}

\section{Cobordisms between Clifford and Chekanov tori}
\label{sec:cobord-betw-cliff}

Assume that a Lagrangian $L$ possesses a \emph{Lagrangian} surgery
disc $D$ and let $L'$ be the result of $(n-1)$-antisurgery on $L$
along $D$. As discussed in Section \ref{sec:lagrangian-0-surgery},
this implies that $L$ can be obtained by resolving a double point of
$L'$ and then applying a Lagrangian isotopy. The positive end of the
desingularized antisurgery cobordism $V: L^\natural \leadsto L$ is
also a resolution of the same double point of $L'$ by Lagrangian
0-surgery; as stated in Proposition \ref{prop:distinct_0_surgeries},
$L$ and $L^\natural$ belong to distinct families of such
resolutions.

\subsection{Clifford and Chekanov tori}
\label{sec:cliff-chek-tori}
We now specialize to the case in which $L$ is the Whitney sphere
$S^2_{Wh}$ in $\mathbb C^2$. Resolving its double point by Lagrangian
0-surgery produces a Lagrangian torus which is monotone, since the
boundary of the Lagrangian disc created when performing this surgery
generates one summand of $H_1$ of the torus and has Maslov index 0.

\begin{figure}[b]
    \centering
    \includegraphics[scale=1]{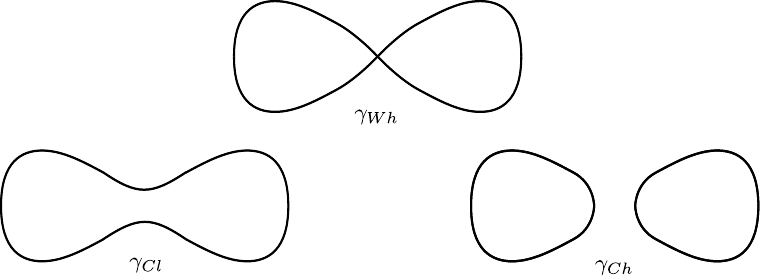}
    \caption{The curves used for constructing $S^2_{Wh}$, $T^2_{Cl}$
        and $T^2_{Ch}$.}
    \label{fig:rotation-curves}
\end{figure}

In fact, the two topologically different types of resolving the double
point yield a Clifford type torus $T_{Cl}^2$ in one case and a
Chekanov type torus $T_{Ch}^2$ in the other case. One can see that by
viewing each of the three Lagrangians as obtained by rotating certain
curves $\gamma: S^1 \to \mathbb C$, i.e., as
\begin{equation*}
    L_\gamma = \{(\gamma(e^{is}) e^{it}, \gamma(e^{is}) e^{-it}) \in \mathbb
    C^2~|~ s,t \in [0,2\pi]\}.
\end{equation*}
To obtain $S^2_{Wh}$ in that way, one uses a figure-8 curve
$\gamma_{Wh}$ with a double point at the origin and symmetric with
respect to the involution $z \mapsto -z$ (to be precise, this yields
the image of the standard Whitney immersion
$S^2 \to T^*\mathbb R^2 \cong \mathbb C^2$ given in
\eqref{eq:Whitney-immersion} under a linear symplectomorphism).
Resolving the double point of $S^2_{Wh}$ by 0-surgery has the same
result as resolving the double point of the figure-8 curve in such a
way that it stays symmetric with respect to $z \mapsto -z$, and then
rotating the resulting curve. The two different ways of performing
this surgery yield a connected curve $\gamma_{Cl}$ enclosing the
origin in one case, and a disconnected curve $\gamma_{Ch}$ whose
components do not enclose the origin in the other case, see Figure
\ref{fig:rotation-curves}. The corresponding Lagrangian tori are
$T^2_{Cl}$ resp.\ $T^2_{Ch}$ up to Hamiltonian isotopy, see e.g.\
\cite{eliashberg-polterovich1997, gadbled_exotic_2013}.

For better compatibility of this description with our setting, it is
useful to note that one can also view the Lagrangians
$S^2_{Wh}, T^2_{Cl}$ and $T^2_{Ch}$ as orbits of curves
$\gamma \subset T^*\mathbb R_1$ under the standard $SO(2)$-action on
$T^*\mathbb R^2$ given by $A \cdot (\x, \y) = (A\x, A\y)$ for
$A \in SO(2)$. Namely, there exists a linear symplectomorphism
$\mathbb C^2 \to T^*\mathbb R^2$ which is equivariant with respect to
the $S^1$-action on $\mathbb C^2$ given by
$e^{it} \cdot (z_1, z_2) = (e^{it} z_1, e^{-it} z_2)$ and the standard
$SO(2)$-action on $T^*\mathbb R^2$ (identifying $SO(2) \cong S^1$ in
the usual fashion), and that takes the plane
$\{ (z, z) \in \mathbb C^2 ~|~ z \in \mathbb C\}$ to the
$(x_1, y_1)$-plane $T^* \mathbb R_1 \subset T^*\mathbb
R^2$. Explicitly, this symplectomorphism is given by
\begin{equation*}
    \big(a_1 + i b_1, a_2 + i b_2\big) \mapsto \frac{1}{2} \big(a_1 + a_2, b_1 -
    b_2, b_1 + b_2, -a_1 + a_2 \big)
\end{equation*}
for $(z_1 = a_1 + ib_1, z_2 = a_2 + ib_2) \in \mathbb C^2$; on the
right-hand side, the first two coordinates correspond to the
0-section, and the last two to the fibre.

\subsection{Cobordisms between the tori}
\label{sec:cobordisms-clifford-chekanov}
In what follows, we denote by $T_{Cl}^2(A)$ and $T_{Ch}^2(A)$ Clifford
and Chekanov tori for which a disc of Maslov index 2 has area $A > 0$
(i.e., for which the monotonicity constant is $\frac{1}{2}A$);
moreover, we denote by $S^2_{Wh}(A)$ the Whitney sphere for which a
generator of $H_2(T^*\mathbb R^2,S^2_{Wh}) \cong \mathbb Z$ has area
$A$. In all cases, $A$ is half of the area bounded by the respective
curves in Figure \ref{fig:rotation-curves}.

The following theorem constructs cobordisms between a Chekanov and a
Clifford torus by first applying 1-antisurgery to the Chekanov torus
such as to obtain a Whitney sphere, and then desingularizing the
corresponding antisurgery cobordism to get a cobordism whose positive
end is a Clifford torus.

\begin{figure}[t]
    \centering
    \includegraphics[scale=1.00]{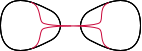}
    \caption{Antisurgering the Chekanov torus}
    \label{fig:surgery-tori}
\end{figure}

\begin{thm}
    \label{thm:Clifford-Chekanov}
    For every choice of $a < A$ with $a/A$ sufficiently close to 1,
    there exists a Lagrangian cobordism
    $T^2_{Cl}(a) \leadsto T^2_{Ch}(A)$ which as smooth a manifold is
    obtained from $[0,1] \times T^2$ by successively attaching a
    2-handle and a 1-handle.
\end{thm}
\begin{proof}
    Consider a curve $\gamma_{Ch}$ of the type required for
    constructing a Chekanov torus $T^2_{Ch}(A)$, i.e., as in the lower
    right part of Figure \ref{fig:rotation-curves} and such that the
    area bounded by each component is $A$. By performing 0-antisurgery
    on $\gamma_{Ch}$, we obtain an immersed curve $\gamma_{Wh}$ which
    after rotation yields a Whitney sphere $S^2_{Wh}(\alpha)$, see
    Figure \ref{fig:surgery-tori}. It follows from the computation in
    Section \ref{sec:lagr-vs.-hamilt} that the areas $A$ and $\alpha$
    are related by $\alpha = A - 2 \varepsilon^{3/2}$, where
    $\varepsilon$ is the size parameter of the antisurgery model we
    implant. If we subsequently resolve the double point of
    $\gamma_{Wh}$ by implanting a local model for Lagrangian 0-surgery
    of size $\delta$ (in the sense of Definition
    \ref{defn:size-of-0-surgery}), we obtain a curve $\gamma_{Cl}$
    bounding area $2(\alpha + \delta)$, which after rotating yields a
    Clifford torus $T^2_{Cl}(\alpha + \delta)$.

    By inspecting the definitions of the antisurgery models in Section
    \ref{sec:immersed-cobordism}, one sees easily that the model for
    1-antisurgery in $T^*\mathbb R^2$ is the $SO(2)$-orbit of the
    model for $0$-antisurgery in $T^*\mathbb R$, viewing the latter as
    living in $T^*\mathbb R_1 \subset T^*\mathbb R^2$. Similarly, the
    Lagrangians $T^2_{Ch}(A)$ and $S^2_{Wh}(\alpha)$ are
    $SO(2)$-orbits of the curves $\gamma_{Ch}, \gamma_{Wh}$, as noted
    at the end of Section \ref{sec:cliff-chek-tori}. Combining these
    statements and the fact that $\gamma_{Wh}$ is the result of
    0-antisurgery on $\gamma_{Ch}$, we see that $S^2_{Wh}(\alpha)$ is
    the result of 1-antisurgery on $T^2_{Ch}(A)$ (and thus
    $T^2_{Ch}(A)$ is the result of 0-surgery on
    $S_{Wh}^2(\alpha)$). 

    Consider now the desingularized antisurgery cobordism correponding
    to the 1-antisurgery that takes $T^2_{Ch}(A)$ to
    $S^2_{Wh}(\alpha)$. This cobordism has as its negative end
    $T_{Ch}^2(A)$ and as its positive end a resolution by 0-surgery
    of $S^2_{Wh}(\alpha)$ of the form $T^2_{Cl}(\alpha + \delta)$ for
    small $\delta > 0$ (using Proposition
    \ref{prop:distinct_0_surgeries}), i.e., it is a cobordism
    $T^2_{Cl}(\alpha + \delta) \leadsto T^2_{Ch}(A)$. Moreover, it has
    the claimed topology by Theorem
    \ref{main-thm:embedded-cob}. Recalling the relation
    $\alpha = A - 2 \varepsilon^{3/2}$, one sees that with this
    construction one can obtain a cobordism
    $T^2_{Cl}(a) \leadsto T^2_{Ch}(A)$ for any $0 < a < A$ for which
    $a/A$ is sufficiently close to 1 by making $\varepsilon$ and
    $\delta$ sufficiently small.
\end{proof}

It seems likely that one can deform the curve $\gamma_{Ch}$ and the
antisurgery model in such a way that the area parameter $\alpha$
describing the size of the Whitney sphere $S^2_{Wh}(\alpha)$ that
comes up in the proof of Theorem \ref{thm:Clifford-Chekanov} is
arbitrarily close to 0. This would imply the existence of a cobordism
$T^2_{Cl}(a) \leadsto T^2_{Cl}(A)$ for \emph{any} choice of
$0 < a < A$ (cf.\ the discussion in Section
\ref{sec:monot-outp-surg}).
    
On the other hand, it is \emph{not} possible to construct a cobordism
$T^2_{Cl}(A) \leadsto T^2_{Ch}(A)$, i.e., between Clifford and
Chekanov type tori of the same size, by our method. Indeed, such a
cobordism would be \emph{monotone} by Proposition
\ref{prop:H_2-surjectivity} below, which would imply equality of
counts of pseudoholomorphic discs of Maslov index 2 through a given
point on both ends, as first observed in \cite{Chekanov--Lagr-cobs}
(see also \cite{Biran-Cornea--Lag-Cob-I, Biran-Cornea--Lag-Cob-II}).
However, it is well known that these counts are different for
$T_{Cl}^2$ and $T_{Ch}^2$. The same argument shows that one cannot
build a Lagrangian cobordism between the \emph{monotone} Clifford and
Chekanov tori in $\mathbb CP^2$ or $S^2 \times S^2$ by this method,
since the monotonicity constant of any monotone Lagrangian there is
determined by that of the ambient manifold (in particular, it is the
same for Clifford and Chekanov type tori).

\begin{prop}
    \label{prop:H_2-surjectivity}
    Let $V^\natural: L^\natural \leadsto L$ be a Lagrangian cobordism
    obtained by desingularizing a cobordism arising from
    $(n-1)$-antisurgery on an embedded Lagrangian submanifold $L$ of a
    symplectic manifold $(M, \omega)$. If $H_1(M) = 0$, then the map
    $H_2(M,L) \oplus H_2(M,L^\natural) \to H_2(T^*\mathbb R \times M,
    V^\natural)$ induced by the inclusions of the ends is surjective.
    In particular, if $L$ and $L^\natural$ are both monotone with the
    same monotonicity constant, then $V^\natural$ is monotone as well.
\end{prop} 
\begin{proof}
    As a smooth manifold, $V^\natural$ is obtained from the cylinder
    $[0, 1] \times L$ by successively attaching an $n$-handle and a
    1-handle. One can see from this description that the map
    $H_1(L) \oplus H_1(L^\natural) \to H_1(V^\natural)$ induced by the
    inclusion of the ends is surjective; in fact, there exist a
    generator $\gamma \in H_1(L^\natural)$ such that the restriction
    of this map to
    $H_1(L) \oplus \mathbb Z \gamma \to H_1(V^\natural)$ is
    surjective. Consider now the following commutative diagram, where
    the horizontal arrows come from the exact sequences of the various
    pairs and the vertical ones are induced by inclusions:
    \begin{equation*}
        \begin{aligned}
            \xymatrix@C=1.65em{H_2(M)^{\oplus 2}\ar[r]\ar[d]&H_2(M,L) \oplus
                H_2(M,L^\natural) \ar[r]\ar[d]& H_1(L)\oplus
                H_1(L^\natural) \ar[r]\ar[d]&
                H_1(M)^{\oplus 2} \ar[d]\\
                H_2(T^*\mathbb R \times M) \ar[r] & H_2(T^*\mathbb R \times
                M, V^\natural) \ar[r] & H_1(V^\natural) \ar[r] &
                H_1(T^*\mathbb R \times M)}
        \end{aligned}
    \end{equation*}
    Since the first and third vertical maps are surjective and the
    fourth is an isomorphism as $H_1(M) = 0$ by assumption, it follows
    from the 4-lemma that the map
    $H_2(M,L) \oplus H_2(M,L^\natural) \to H_2(T^*\mathbb R \times M,
    V^\natural)$ is also surjective.

    The last statement follows easily by observing that the pullbacks
    of the area and Maslov homomorphisms
    $H_2(T^*\mathbb R \times M, V^\natural) \to \mathbb R$ by the
    natural inclusions $L, L^\natural \hookrightarrow V^\natural$ of
    the ends are the corresponding homomorphisms
    $H_2(M, L) \to \mathbb R$ resp.\
    $H_2(M, L^\natural) \to \mathbb R$.
\end{proof}

\subsection{Successive antisurgery/surgery for Clifford and Chekanov tori}
\label{sec:surg-cliff-chek}
As discussed in Sections \ref{sec:lagrangian-0-surgery-1} and
\ref{sec:lagr-vs.-hamilt}, the result of successively applying
$(n-1)$-antisurgery to a Lagrangian $L$ and then $0$-surgery to the
resulting $L'$ leads to a Lagrangian $L^\natural$ which is Lagrangian
isotopic to the original $L$ (for one of the two topologically
different ways of implanting the local model for 0-surgery), but not
Hamiltonian isotopic to $L$. In particular, if $L$ is a torus of
Clifford or Chekanov type, then the resulting $L^\natural$ obtained
that way is again a torus of the same type. Assume that we perform the
two operations as described in Section
\ref{sec:cobordisms-clifford-chekanov}, i.e., by first operating on
the level of curves and then rotating. Note that if we start with a
Chekanov torus $T^2_{Cl}(A)$, then the Whitney sphere obtained by
antisurgery has area less than $A$, as can be seen from Figure
\ref{fig:surgery-tori}, and hence the resulting $T^2_{Ch}(a)$ has a
\emph{smaller} area parameter than the original one.  On the other
hand, if we start with a Clifford torus $T^2_{Cl}(a)$, then the
resulting $T^2_{Cl}(A)$ has \emph{larger} area parameter. It would be
interesting to investigate if this observation can be used to give
another proof of the fact that $T_{Cl}^2(A)$ and $T_{Ch}^2(A)$ are not
Hamiltonian isotopic.

\section*{Acknowledgements}
\label{sec:acknowledgements}
I am particularly grateful to Paul Biran for his continued interest in
this project and for providing plenty of feedback, as well as Mads
Bisgaard and Jeff Hicks for very useful discussions which helped me
spot a mistake in an earlier version of this article. I also would
like to thank Denis Auroux, Fran\c{c}ois Charette, Octav Cornea,
Georgios Dimitroglou Rizell, Tobias Ekholm, Jonny Evans, Leonid
Polterovich, Paul Seidel, Dmitry Tonkonog and Weiwei Wu for related
discussions and comments. Finally, I am grateful for the thoughtful
comments of an anonymous referee that have helped to improve the
article.

This article was mainly written while I was employed as a CIRGET
postdoctoral fellow at the Centre de Recherches Mathématiques in
Montreal, and partly during a stay at the Institut Mittag-Leffler in
Stockholm. I thank both institutions for their hospitality.

\bibliographystyle{alpha}
\bibliography{antisurgery}

\end{document}